\input amssym.def
\input amssym.tex

\def\item#1{\vskip1.3pt\hang\textindent {\rm #1}}


\newskip\litemindent
\litemindent=0.7cm  
\def\Litem#1#2{\par\noindent\hangindent#1\litemindent
\hbox to #1\litemindent{\hfill\hbox to \litemindent
{\ninerm #2 \hfill}}\ignorespaces}
\def\litem{\Litem1}

\tolerance=300
\pretolerance=200
\hfuzz=1pt
\vfuzz=1pt

\hoffset=0in
\voffset=0.5in

\hsize=5.8 true in
\vsize=9.2 true in
\parindent=25pt
\mathsurround=1pt
\parskip=1pt plus .25pt minus .25pt
\normallineskiplimit=.99pt

\countdef\revised=100
\mathchardef\emptyset="001F 
\chardef\ss="19
\def\3{\ss}
\def\anf{$\lower1.2ex\hbox{"}$}
\def\frac#1#2{{#1 \over #2}}
\def\>{>\!\!>}
\def\<{<\!\!<}

\def\into{\hookrightarrow}

\def\ssssarr{\hbox to 15pt{\rightarrowfill}}
\def\sssarr{\hbox to 20pt{\rightarrowfill}}
\def\ssarr{\hbox to 30pt{\rightarrowfill}}
\def\sarr{\hbox to 40pt{\rightarrowfill}}
\def\arr{\hbox to 60pt{\rightarrowfill}}
\def\larr{\hbox to 60pt{\leftarrowfill}}
\def\Arr{\hbox to 80pt{\rightarrowfill}}

\def\ad{\mathop{\rm ad}\nolimits}

\def\Ad{\mathop{\rm Ad}\nolimits}

\def\algint{\mathop{\rm algint}\nolimits}
\def\Aut{\mathop{\rm Aut}\nolimits}

\def\PSU{\mathop{\rm PSU}\nolimits}

\def\Fix{\mathop{\rm Fix}\nolimits}
\def\End{\mathop{\rm End}\nolimits}
\def\Ext{\mathop{\rm Ext}\nolimits}

\def\GL{\mathop{\rm GL}\nolimits}
\def\Gr{\mathop{\rm Gr}\nolimits}

\def\Herm{\mathop{\rm Herm}\nolimits}
\def\Hom{\mathop{\rm Hom}\nolimits}%
\def\id{\mathop{\rm id}\nolimits} 
\def\im{\mathop{\rm im}\nolimits}

\def\Im{\mathop{\rm Im}\nolimits}

\def\Int{\mathop{\rm int}\nolimits}

\def\OO{\mathop{\rm O{}}\nolimits}

\def\PSL{\mathop{\rm PSL}\nolimits}
\def\PSp{\mathop{\rm PSp}\nolimits}
\def\PSO{\mathop{\rm PSO}\nolimits}
\def\PU{\mathop{\rm PU}\nolimits}

\def\rank{\mathop{\rm rank}\nolimits}
\def\rk{\mathop{\rm rank}\nolimits}

\def\Re{\mathop{\rm Re}\nolimits}
\def\rk{\mathop{\rm rk}\nolimits}

\def\SL{\mathop{\rm SL}\nolimits}
\def\SO{\mathop{\rm SO}\nolimits}
\def\span{\mathop{\rm span}\nolimits}

\def\Sym{\mathop{\rm Sym}\nolimits}

\def\Sp{\mathop{\rm Sp}\nolimits}

\def\Spin{\mathop{\rm Spin}\nolimits}

\def\SU{\mathop{\rm SU}\nolimits}
\def\sup{\mathop{\rm sup}\nolimits}

\def\tr{\mathop{\rm tr}\nolimits}

\def\UU{\mathop{\rm U{}}\nolimits}

\def\0{{\bf 0}}
\def\1{{\bf 1}}

\def\g{{\frak g}}
\def\gl{{\frak {gl}}}

\def\su{{\frak {su}}}

\def\sL{{\frak {sl}}}

\def\uu{{\frak u}}

\def\L{\mathop{\bf L{}}\nolimits}

\def\C{{{\Bbb C}{\mskip+1mu}}} 
\def\K{{{\Bbb K}{\mskip+2mu}}} 
\def\H{{{\Bbb H}{\mskip+2mu}}} 

\def\R{{\Bbb R}}
\def\Z{{\Bbb Z}}
\def\N{{\Bbb N}}

\def\E{{\Bbb E}}
\def\F{{\Bbb F}}
\def\K{{\Bbb K}}

\def\P{{\Bbb P}}

\def\SS{{\Bbb S}}
\def\T{{\Bbb T}}

\def\:{\colon}  
\def\.{{\cdot}}
\def\|{\Vert}
\def\bsk{\bigskip}

\def\giantskip{\vskip2\bigskipamount}
\def\gsk{\giantskip}
\def \la {\langle}

\def\msk{\medskip}
\def \ra {\rangle}
\def \res {\!\mid\!\!}

\def\bbr{\bigbreak}
\def\giantbreak{\par \ifdim\lastskip<2\bigskipamount \removelastskip
          \penalty-400 \giantskip\fi}

\def\nin{\noindent}
\def\cen{\centerline}
\def\pagebreak{\vskip 0pt plus 0.0001fil\break}
\def\linebreak{\break}

\def\eps{\varepsilon}
\def\epsilon{\varepsilon}

\def\nin{\noindent}
\def\oline{\overline}

\def\pder#1,#2,#3 { {\partial #1 \over \partial #2}(#3)}
\def\pde#1,#2 { {\partial #1 \over \partial #2}}
\def\phi{\varphi}


\def\subeq{\subseteq}
\def\supeq{\supseteq}

\def\tilde{\widetilde}

\font\ninerm=cmr9
\font\eightrm=cmr8

\font\eightbf=cmbx8


\font\smc=cmcsc10
\font\bfone=cmbx10 scaled\magstep1 
\font\bftwo=cmbx10 scaled\magstep2 

\def\qed{{\unskip\nobreak\hfil\penalty50\hskip .001pt \hbox{}\nobreak\hfil
           \vrule height 1.2ex width 1.1ex depth -.1ex
            \parfillskip=0pt\finalhyphendemerits=0\medbreak}\rm}

\def\qeddis{\eqno{\vrule height 1.2ex width 1.1ex depth -.1ex} $$
                    \medbreak\rm}

\def\Lemma #1. {\bigbreak\vskip-\parskip\noindent{\bf Lemma #1.}\quad\it}

\def\Sublemma #1. {\bigbreak\vskip-\parskip\noindent{\bf Sublemma #1.}\quad\it}

\def\Proposition #1. {\bigbreak\vskip-\parskip\noindent{\bf Proposition #1.}
\quad\it}

\def\Corollary #1. {\bigbreak\vskip-\parskip\nin{\bf Corollary #1.}
\quad\it}

\def\Theorem #1. {\bigbreak\vskip-\parskip\noindent{\bf Theorem #1.}
\quad\it}

\def\Definition #1. {\rm\bigbreak\vskip-\parskip\noindent
{\bf Definition #1.}
\quad}

\def\Remark #1. {\rm\bigbreak\vskip-\parskip\noindent{\bf Remark #1.}\quad}

\def\Example #1. {\rm\bigbreak\vskip-\parskip\noindent{\bf Example #1.}\quad}
\def\Examples #1. {\rm\bigbreak\vskip-\parskip\noindent{\bf Examples #1.}\quad}

\def\Problems #1. {\bigbreak\vskip-\parskip\noindent{\bf Problems #1.}\quad}
\def\Problem #1. {\bigbreak\vskip-\parskip\noindent{\bf Problem #1.}\quad}
\def\Exercise #1. {\bigbreak\vskip-\parskip\noindent{\bf Exercise #1.}\quad}

\def\Conjecture #1. {\bigbreak\vskip-\parskip\noindent{\bf Conjecture 
#1.}\quad}

\def\Proof#1.{\rm\par\ifdim\lastskip<\bigskipamount\removelastskip\fi\smallskip
             \noindent {\bf Proof.}\quad}

\def\Axiom #1. {\bigbreak\vskip-\parskip\noindent{\bf Axiom #1.}\quad\it}

\def\Satz #1. {\bigbreak\vskip-\parskip\noindent{\bf Satz #1.}\quad\it}

\def\Korollar #1. {\bbr\vskip-\parskip\nin{\bf Korollar #1.} \quad\it}

\def\Folgerung #1. {\bbr\vskip-\parskip\nin{\bf Folgerung #1.} \quad\it}

\def\Folgerungen #1. {\bbr\vskip-\parskip\nin{\bf Folgerungen #1.} \quad\it}

\def\Bemerkung #1. {\rm\bigbreak\vskip-\parskip\noindent{\bf Bemerkung #1.}
\quad}

\def\Beispiel #1. {\rm\bigbreak\vskip-\parskip\noindent{\bf Beispiel #1.}\quad}
\def\Beispiele #1. {\rm\bigbreak\vskip-\parskip\noindent{\bf 
Beispiele #1.}\quad}
\def\Aufgabe #1. {\rm\bigbreak\vskip-\parskip\noindent{\bf Aufgabe #1.}\quad}
\def\Aufgaben #1. {\rm\bigbreak\vskip-\parskip\noindent{\bf Aufgabe #1.}\quad}

\def\Beweis#1. {\rm\par\ifdim\lastskip<\bigskipamount\removelastskip\fi
            \smallskip\noindent {\bf Beweis.}\quad}

\nopagenumbers

\def\date{\ifcase\month\or January\or February \or March\or April\or May
\or June\or July\or August\or September\or October\or November
\or December\fi\space\number\day, \number\year}

\def\title{Title ??}
\def\author{Author ??}

\def\thanks#1{\footnote*{\eightrm#1}}

\def\rightheadline{\hfil{\eightrm\title}\hfil\tenbf\folio}
\def\leftheadline{\tenbf\folio\hfil{\eightrm\author}\hfil}
\headline={\vbox{\line{\ifodd\pageno\rightheadline\else\leftheadline\fi}}}

\def\firstheadline{}
\def\firstfootline{\cen{\rm\folio}}

\def\seite #1 {\pageno #1
                \headline={\ifnum\pageno=#1 \firstheadline
                \else\ifodd\pageno\rightheadline\else\leftheadline\fi\fi}
                \footline={\ifnum\pageno=#1 \firstfootline\else{}\fi}}

\newdimen\dimenone
  \def\checkleftspace#1#2#3#4{
  \dimenone=\pagetotal
  \advance\dimenone by -\pageshrink   
  \ifdim\dimenone>\pagegoal          
    \else\dimenone=\pagetotal
         \advance\dimenone by \pagestretch
         \ifdim\dimenone<\pagegoal
           \dimenone=\pagetotal
           \advance\dimenone by#1         
           \setbox0=\vbox{#2\parskip=0pt                
                      \hyphenpenalty=10000
                      \rightskip=0pt plus 5em
                      \noindent#3 \vskip#4}    
         \advance\dimenone by\ht0
         \advance\dimenone by 3\baselineskip
         \ifdim\dimenone>\pagegoal\vfill\eject\fi
           \else\eject\fi\fi}


\def\subheadline #1{\nin\bigbreak\vskip-\lastskip
       \checkleftspace{0.9cm}{\bf}{#1}{\medskipamount}
           \indent\vskip0.7cm\centerline{\bf #1}\medskip}
\def\subsection{\subheadline}

\def\lsubheadline #1 #2{\bigbreak\vskip-\lastskip
       \checkleftspace{0.9cm}{\bf}{#1}{\bigskipamount}
          \vbox{\vskip0.7cm}\cen{\bf #1}\msk \cen{\bf #2}\bsk}

\def\sectionheadline #1{\bigbreak\vskip-\lastskip
       \checkleftspace{1.1cm}{\bf}{#1}{\bigskipamount}
          \vbox{\vskip1.1cm}\cen{\bfone #1}\bsk}
\def\section{\sectionheadline}

\def\lsectionheadline #1 #2{\bigbreak\vskip-\lastskip
       \checkleftspace{1.1cm}{\bf}{#1}{\bigskipamount}
          \vbox{\vskip1.1cm}\cen{\bfone #1}\msk \cen{\bfone #2}\bsk}

\def\lchapterheadline #1 #2{\bigbreak\vskip-\lastskip\indent\vskip3cm
                        \cen{\bftwo #1} \msk \cen{\bftwo #2} \gsk}
\def\llsectionheadline #1 #2 #3{\bigbreak\vskip-\lastskip\indent\vskip1.8cm
\cen{\bfone #1} \msk \cen{\bfone #2} \msk \cen{\bfone #3} \nobreak\bsk\nobreak}


\newtoks\literat
\def\[#1 #2\par{\literat={#2\unskip.}%
\hbox{\vtop{\hsize=.15\hsize\nin [#1]\hfill}
\vtop{\hsize=.82\hsize\nin\the\literat}}\par
\vskip.3\baselineskip}

\def\references{
\sectionheadline{\bf References}
\frenchspacing

\entries\par}

\mathchardef\emptyset="001F
\def\address{Author: \tt$\backslash$def$\backslash$address$\{$??$\}$}

\def\abstract #1{{\narrower\baselineskip=10pt{\noindent
\eightbf Abstract.\quad \eightrm #1 }
\bigskip}}

\def\firstpage{\nin
{\obeylines \parindent 0pt }
\vskip2cm
\centerline{\bfone\title}
\gsk
\centerline{\bf\author}
\vskip1.5cm \rm}

\def\addresstwo{}

\def\dlastpage{\par\vbox{\vskip1cm\nin
\line{
\vtop{\hsize=.5\hsize{\parindent=0pt\baselineskip=10pt\nin\address}}
\quad
\vtop{\hsize=.42\hsize\nin{\parindent=0pt
\baselineskip=10pt\addresstwo}}
\hfill} }}

\def\Box #1 { \msk\par\nin
\centerline{
\vbox{\offinterlineskip
\hrule
\hbox{\vrule\strut\hskip1ex\hfil{\smc#1}\hfill\hskip1ex}
\hrule}\vrule}\msk }

\def\adots{\mathinner{\mkern1mu\raise1pt\vbox{\kern7pt\hbox{.}}
                         \mkern2mu\raise4pt\hbox{.}
                         \mkern2mu\raise7pt\hbox{.}\mkern1mu}}


\def\title{Orbits of triples in the Shilov boundary of a bounded symmetric
domain}
\def\author{Jean-Louis Clerc\footnote{*}{\rm
The first author acknowledges partial support from the EU (TMR
Network Harmonic Analysis and Related Problems)}, Karl-Hermann Neeb
}
\def\date{September 29, 2005}
\def\rightheadline{\tenbf\folio\hfil{\tt triples.tex}\hfil\eightrm\date}
\def\leftheadline{\tenbf\folio\hfil{\rm\title}\hfil\eightrm\date}

\def\address
{Jean-Louis Clerc

Institut Elie Cartan

Facult\'e des Sciences, Universit\'e Nancy I

B.P. 239

F - 54506 Vand\oe uvre-l\`es-Nancy Cedex

France

clerc@iecn.u-nancy.fr

}

\def\addresstwo
{Karl-Hermann Neeb

Technische Universit\"at Darmstadt

Schlossgartenstrasse 7

D-64289 Darmstadt

Germany

neeb@mathematik.tu-darmstadt.de}

\def\Face{\mathop{\rm Face}\nolimits}

\def\mod{\mathop{\rm mod}\nolimits}
\def\half{\textstyle{1\over 2}}

\firstpage

\abstract{ Let ${\cal D}$ be a bounded symmetric domain of tube type, $S$ its
Shilov boundary, and $G$ the neutral component of
its group of biholomorphic transforms. We classify
the orbits of $G$ in the set $S\times S\times S$.\hfill\break
Keywords: bounded symmetric domain, tube type domain, Shilov boundary,
face, Maslov index, flag manifold, Jordan triples \hfill\break
MSC: 32M15, 53D12, 22F30\hfill}

\subheadline{Introduction}

Let $\cal D$ be a bounded symmetric domain, realized as a circular
domain in a (finite dimensional) complex vector
space $V$. Let
$G := \Aut({\cal D})_0$ be the identity component of its group of
biholomorphic transforms
of ${\cal D}$ and let
$S$ be the Shilov boundary of ${\cal D}$.
The action of any element of $G$ extends to a neigbourhood of
$\overline{\cal D}$, and hence $G$ acts on $S$. It is well known that this
action is transitive.
The main result of the present paper is a classification of the $G$-orbits in
the set $S\times S\times S$ of triples in $S$, when
$\cal D$ is of {\it tube type\/}.

The action of $G$ on $S\times S$ can be easily studied as an
application of Bruhat theory,
and the description of the orbits is the same, whether $\cal D$ is of
tube type or not. But
for triples, there is a drastic difference between tube type domains
and non tube type
domains. In the first case, there is a finite number of orbits in
$S\times S\times
S$, whereas there are an infinite number of orbits for a non tube type domain.

Let $r$ be the rank of $\cal D$. The notion of {\it $r$-polydisc\/}  (and its
corresponding Shilov boundary called
{\it $r$-torus\/}) plays an important role in the analysis of the orbits. On
one hand they
are the ``complexifications'' of the {\it maximal flats\/} of ${\cal
D}$ (in the sense of the
geometry of Riemannian symmetric spaces). On the other hand, a $r$-polydisc in
the usual sense is a set of the form
$$\Delta^r = \Big\{\sum_{j=1}^r \zeta_jx_j\: |\zeta_j|<1 , 1\leq
j\leq r\Big\}\ ,$$
where the $x_j$ are linearly independent elements in $V$. The space
$V$ has a natural
structure of a positive
hermitian Jordan triple system, and in particular, it has a natural (Banach)
norm, called the {\it spectral norm\/},
for which the domain $\cal D$ is realized as the open unit ball.
One of the results we prove is that such a polydisc, constructed on vectors
$x_j$ of norm $1$  lies in $\cal D$ if and only if the $(x_j)_{1\leq
j\leq r}$ form a
{\it Jordan frame\/} for $V$.

Fix an $r$-torus $T \subeq S$ arising as the Shilov boundary of an
$r$-polydisc associated to a Jordan
frame. The main step towards the classification of the orbits of $G$ in
$S\times S\times S$ is the result that any triple in
$S$ can be sent by an element of
$G$ to a triple in $T$. This requires that $\cal D$ is of tube type,
and this property
really distinguishes tube type domains from non tube type domains.
Once this result is
obtained, the classification becomes easy, because the problem is
reduced to the case of a
polydisc, and further, using the product structure, to the case of
the unit disc in
$\Bbb C$, where the situation is easy to analyze. The generalized Maslov index
(see [C\O{}01], [Cl04]) comes in as a subtle invariant for triples.

A special case of this theorem was known before. If $\cal D$ is the
Siegel domain (the unit ball in the space of complex symmetric matrices
$\Sym_r(\Bbb C)$), then the group
$G$ is the projective symplectic group $\PSp_{2r}(\R) :=
\Sp_{2r}(\R)/\{\pm\1\}$,
and the Shilov boundary of $\cal D$ can be
identified with the Lagrangian manifold (the set of Lagrangian subspaces
of $\Bbb R^{2r}$). Then the orbits of triples of
Lagrangians have been described (see [KS90, p.492]), using linear
symplectic algebra
techniques. Related results can be found in
[FMS04],  and in particular their
Proposition 4.3 (which they deduce from [KS90]) is,
for this specific example, equivalent to our Theorem III.1.
The main point of [FMS04] is a description of the orbits of the action
of the maximal compact subgroup group $\UU_n(\C)$ of $\Sp_{2n}(\R)$
on triples of
Lagrangians are classified, but this is a different problem.

As explained in the appendix, the bounded symmetric domain of tube 
type  can be described
in terms of euclidean Jordan algebras $E$. More precisely, the irreducible ones
are in one-to-one correspondence with simple euclidean
Jordan algebras. From the table in [FK94, p.~213]
(see also [Be00]) it is easy to give the following table, where for each simple
euclidean Jordan algebra $E$, we list  the group $L$ of linear
transforms of $E$ preserving the cone
$\Omega$, the group $G$ of holomorphic diffeomorphisms of the bounded 
symmetric domain
$\cal D$ and the Shilov boundary $S$ as compact Riemannian symmetric 
space. There are four
infinite series and one exceptional case. From the point of view of 
flag manifold (see
below), $S$ is realized as $G/P$, where the (maximal) parabolic 
subgroup $P$ is the
semi-direct product of $L$ (Levi component) and $E$ (unipotent radical).
\msk
\centerline {Table 1.}
\msk
\def\tablerule{\noalign{\hrule}}
\vbox{\tabskip=0pt\offinterlineskip
\halign to 14.9cm{
\strut#&
\vrule# \tabskip=0em plus 2em&#\hfil&
\vrule\
\vrule#&\hfil#\hfil&
\vrule#&\hfil#\hfil&
\vrule#&\hfil#\hfil&
\vrule#&\hfil#\hfil&
\vrule#&\hfil#\hfil&
\vrule#\tabskip=0pt\cr\tablerule

&&\omit\hidewidth $E$ \hidewidth
&&\omit\hidewidth $\Sym_n(\R)$  \hidewidth
&&\omit\hidewidth $\Herm_n(\C)$ \hidewidth
&&\omit\hidewidth $\Herm_n(\H)$ \hidewidth
&&\omit\hidewidth $\R^{1,n-1}$ \hidewidth
&&\omit\hidewidth $\Herm_3(\Bbb O)$ \hidewidth
&\cr\tablerule

&& $L$ && $\GL_n(\R)$ && $\GL_n(\C)$ && $\GL_n(\H)$ &&
$\SO_0(1,n-1)\times
\R^* $&& $\E_{6(-26)}\times \Bbb R^*$ &\cr\tablerule

&& $G$ && $\PSp_{2n}(\R)$ && $\PU_{n,n}(\C)$ && $\PSO^*(4n)$
&& $\SO_{2,n}(\R)_0$ && $\E_{7(-25)}$ &\cr\tablerule

&& $S$ && $Lag(\R^{2n})\simeq \UU_n(\C)/\OO_n(\R)$ && $\UU_n(\C)$ &&
$\UU_{2n}(\C)/\SU(n,\H)$  && $\SS^1 \times \SS^{n-1}/\Z_2$ && $\UU(1)\E_6/\F_4$
&\cr\tablerule

\cr}}
\msk

The Shilov boundary $S$ of a bounded domain is in particular a generalized flag
manifold of
$G$, i.e. of the form $G/P$, where $P$ is a parabolic subgroup of $G$. A nice
description of $P$ is obtained after performing a Cayely transform. 
The domain $\cal D$
is transformed to an unbounded domain ${\cal D}^C$ which is a Siegel 
domain of type II and
the group
$P$ is the group of all affine transformations preserving ${\cal 
D}^C$ (see section 1 for
details). The group
$P$ has some specific properties : it is a maximal parabolic subgroup 
of $G$, conjugate to
its opposite. Moreover, one can show that the domain $\cal D$ is of 
tube type if and only if
the unipotent radical $U$ of $P$  is abelian. A natural question 
arises to  which
extent results similar to the ones obtained in this paper could be 
valid for other
generalized flag manifolds. The natural background for this problem 
is the following. If
$P_1, \ldots, P_k$ are parabolic subgroups of a connected semisimple group
$G'$, then the product manifold
$$ M := G'/P_1 \times \ldots \times G'/P_k $$
is called a multiple flag manifold of finite type if the diagonal action of
$G'$ on $M$ has only finitely many orbits.
For $k = 1$ we always have only one orbit, and for $k = 2$ the
finiteness of the set of
orbits follows from the Bruhat decomposition of $G'$.
For $G' = \GL_n(\K)$ or $G' = \Sp_{2n}(\K)$ and $\K$ is an
algebraically closed field of
characteristic zero, it has been shown in [MWZ99/00] that
finite type implies $k \leq 3$, and for $k = 3$ the triples of
parabolics leading to
multiple flag manifolds of finite type are described and the $G'$-orbits
in these manifolds classified. The main technique to achieve these
classifications
was the representation theory of quivers. In [Li94], Littelmann considers
general simple algebraic groups over $\K$ and describes all multiple
flag manifolds of finite type for $k = 3$ under the assumption that
$P_1$ is a Borel subgroup and $P_2$, $P_3$ are maximal parabolics.
Actually Littelmann considers the condition that $B= P_1$ has a dense orbit in
$G'/P_2 \times G'/P_3$, but the results in [Vi86] show that this
implies the finiteness
of the number of $B$-orbits and hence the finiteness of the number of
$G'$-orbits
in $G'/B \times G'/P_2 \times G'/P_3$.
{}From Littelmann's classification one can easily read off that for a
{}maximal parabolic
$P$ in $G'$ the triple product $(G'/P)^3$ is of finite type
if and only if the unipotent radical $U$ of $P$ is abelian and in two
exceptional situations. If $U$ is abelian, then $P$ is the maximal parabolic
defined by a $3$-grading of $\g' = \L(G')$, so that $G'/P$ is the
conformal completion of a
Jordan triple (cf.\ [BN05] for a discussion of such completions in an
abstract setting). This case was also studied in [RRS92].
The first exceptional case,
where $U$ is not
abelian, corresponds to
$G' = \Sp_{2n}(\K)$, where $G'/P = \P_{2n-1}(\K)$ is the projective
space of $\K^{2n}$,
$U$ is the $(2n-1)$-dimensional Heisenberg group and
the Levi complement is $\Sp_{2n-2}(\K) \times \K^\times$. In the
other exceptional
case $G' = \SO_{2n}(\K)$ and $G'/P$ is the highest weight orbit in the
$2^n$-dimensional
spin representation of the covering group $\tilde G' = \Spin_{2n}(\K)$ of $G'$.
Here $U \cong \Lambda^2(\K^n) \oplus \K^n$ also is a $2$-step nilpotent group
and the Levi complement acts like $\GL_n(\K)$ on this group.
It seems that the positive finiteness results have a good chance to carry over
to the split forms of groups over more general fields and in 
particular to $\K =
\R$,
but for real groups not much seems to be known about multiple flag manifolds of
finite type.

If $M = (G'/P)^3$ is a multiple flag manifold of finite type,
$P$ is conjugate to its opposite, and $P = U \rtimes L$
is a Levi decomposition of $P$, then $L$ is the simultaneous
stabilizer of a pair in
$(G'/P)^2$ with an open orbit, and this implies that the conjugation
action of $L$ on $U$
has only finitely many orbits. A closely related but different problem is the
question  when the conjugation action of $P$ on $U$ has finitely many orbits.
According to a  result of Richardson, every parabolic
$P$ has a dense orbit in its unipotent radical $U$,
but this does not imply  finiteness.
For more specific results on this question we refer to
[RRS92], [PR97] and [HR99].

It is perhaps worthwhile to stress
that the proofs we give are one more occurrence of the interaction
between complex analysis
of a bounded symmetric domains and the geometry of convex sets in the
normed space $V$. The
notions of  extremal points or faces of a convex set do play an
important role in our study.

The contents of the paper is as follows. In Section~I we first recall
several facts on
bounded symmetric domains. Our main sources
are Loos' lecture notes [Lo77] and Satake's book [Sa80]. For results concerning
euclidean Jordan algebras we use [FK94]. The main result of Section~I
is a classification
of the $G$-orbits in the set of quasi-invertible (=transversal) pairs in
$\oline{\cal D}$ (Theorem~I.18). For this classification, there would
be no gain in
assuming that
${\cal D}$ is of tube type, so that the theorem is proved in full generality.
However, for the analysis of
$G$-orbits in
$S\times S\times S$ (assuming $\cal D$ to be of tube type), we only
need the  classification result for transversal pairs
$(x,y)$, where
$x
\in S$ and $y\in \overline {\cal D}$. For this case we give a more
direct shorter
proof (see Lemma~I.20), but we think that the general case might also
be useful in
other situations.

The main tool for the classification of $G$-orbits in $S\times
S\times S$ is the characterization of
the transversality relation on $\oline{\cal D}$ in terms of faces of
the compact convex set
$\oline{\cal D}$: Two elements $x,y \in \oline{\cal D}$ are
transversal if and only if
they are not contained in a proper face of $\oline{\cal D}$
(Theorem~II.12). This characterization is also valid for non tube type domains.
A key concept for the classification is the notion of the rank of a face
$F$ of $\oline{\cal D}$. For an irreducible domain ${\cal D}$ of rank
$r$ it takes
values in the set $\{0,1,\ldots, r\}$ and classifies the $G$-orbits
in the set of
faces of ${\cal D}$. It is normalized in such a way that the rank of
$\oline{\cal D}$
as a face is zero and that the extreme points, i.e., the elements in
the Shilov boundary,
are faces of rank $r$. If $\Face(x_1, \ldots, x_n)$ denotes the face
generated by the subset
$\{x_1,\ldots, x_n\}$ of $\oline{\cal D}$, then the function
$$ \oline{\cal D}^n \to \{0,1\ldots, r\}, \quad
(x_1, \ldots, x_n) \mapsto \rank \Face(x_1,\ldots, x_n) $$
is an invariant for the $G$-action on $\oline{\cal D}^n$.

In these terms,
two elements $x,y \in \oline{\cal D}$ are transversal if and only if
$\rank \Face(x,y)= 0$. In Section~III we use this fact to show that
for a domain ${\cal D}$ of tube
type every triple in $S$ is conjugate to a triple in the Shilov
boundary $T$ of a
maximal polydisc $\Delta^r$ defined by a Jordan frame. This reduces
the classification
of $G$-orbits in $S\times S\times S$ to the description of
intersections of these orbits with
$T^3$. This is fully achieved in Section~V by assigning a
$5$-tuple of integer invariants to each orbit and by showing that
triples with the
same invariant lie in the same orbit. The first four components of
this $5$-tuple are
$$ (\rank\Face(x_1,x_2, x_3), \rank\Face(x_1,x_2),
\rank\Face(x_2,x_3), \rank\Face(x_1, x_3)).$$
The fifth component is defined as
the Maslov index $\iota(x_1, x_2, x_3)$ which is discussed in some
detail in Section~IV.
Note that if $(x_1, x_2, x_3)$ is transversal in the sense that all pairs
$(x_1, x_2)$, $(x_2, x_3)$, $(x_3, x_1)$
are transversal, then the first four components
of the invariant vanish, which implies that the $G$-orbits in the set
of transversal triples
are classified by the Maslov index.

We conclude the paper (Section~VI) with a brief discussion of how the
classification of the
$G$-orbits  in $S \times S$ can be interpreted in terms of the Bruhat
decomposition of $G$. Note that, although $S$ always is a generalized 
flag manifold
of the real group $G$, the unipotent radical of the corresponding parabolic
is abelian if and only if the domain ${\cal D}$ is of tube type.
If this is the case, then [Li94] and [RRS92] imply that the complexification
$G_\C$ acts with finitely many orbits on $(G_\C/P_\C)^3$.
For each $G_\C$-orbit $M \subeq (G_\C/P_\C)^3$ meeting the totally
real submanifold $(G/P)^3$ the intersection $M \cap (G/P)^3$ is totally
real in $M$, hence a real form of $M$, and [BS64, Cor.~6.4] implies
that $G$ has only finitely many orbits in $M \cap (G/P)^3$. Alternatively
one can argue with Whitney's Theorem ([Wh57]) that the set of real points of a
complex variety has only finitely many connected components which coincide
with the $G$-orbits in our case. In view
of this argument, its not the finiteness of the $G$-orbits but their
classification and the relation to the Maslov index that is the main
point of the present paper.

In [RRS92, Th.~1.2(b)] one also finds a classification
of the $G_\C$-orbits in $(G_\C/P_\C)^2$ which turns out to be the same
as in the real case (cf.\ Theorem~VI.1).

A final appendix gives a short presentation of the relation between
positive hermitian Jordan triple systems and bounded symmetric
domains on one hand,
between euclidean Jordan algebras and tube type domains on the other hand. This
appendix is designed for readers not familiar with the language of
Jordan algebra and/or Jordan triple system.

We thank
L.~Kramer and H.~Rubenthaler for comments and references  concerning
multiple flag
manifolds of finite type. We also thank several anonymous editors of this
journal for numerous remarks and for pointing out the reference [RRS92].

\sectionheadline{I. Classification of orbits of transversal pairs in
the boundary}

Let ${\cal D}$ be an irreducible circular bounded symmetric domain, so that
${\cal D}$ is the open unit ball for a norm on a complex vector space
$V$ ([Lo77,
Th.4.1]).  In this section we describe the $G$-orbits in the set of
quasi-invertible  pairs of elements in the closure of ${\cal D}$
(cf.~Theorem~I.18
below). Here we do not have to assume that
${\cal D}$ is of tube type.

\nin {\bf I.1. The associated Jordan triple.} On $V$ we consider the
hermitian Jordan
triple product $\{\cdot, \cdot, \cdot \} \: V^3 \to V$ that is
uniquely determined
by the property that for each $v \in V$ the vector field given by the function
$$ \xi_v \: V \to V, \quad z \mapsto v - \{z,v,z\} $$
generates a one-parameter group of automorphisms of ${\cal D}$
([Lo77, Lemma~4.3]). Note that for each $v \in V$ the map $(z,w)
\mapsto \{z,v,w\}$
is symmetric and complex bilinear, and that, for each $a,b\in V$ the
map $z \mapsto
\{a,z,b\}$ is antilinear.  For $x,y \in V$ we define
$Q(x)$ and $x \square y \in \End(V)$ by
$$ Q(x).y := \{x,y,x\} \quad \hbox{ and } \quad x \square y.z := \{x,y,z\}. $$
The Jordan triple structure on $V$ used by
Loos is $\{x,y,z\}' = 2 \{x,y,z\}$, so that his quadratic representation
is given by $Q'(x,y) = 2\{x,y,z\}$, but since Loos defines
$Q'(x)$ as $\half Q'(x,x)$, we obtain the same operators $Q(x) = Q'(x)$.

\msk

\nin {\bf I.2. Tripotents and Peirce decomposition.} An element $e
\in V$ is called a {\it tripotent} if $e = \{e,e,e\}$.
For a tripotent $e \in V$ let
$V_j := V_j(e)$ denote the $j$-eigenspace of the operator $2 e \square e$.
Then we obtain the corresponding {\it Peirce decomposition of $V$}:
$$ V  = V_0 \oplus V_1 \oplus V_2 $$
([Lo77, Th.~3.13]). Since $e \square e$ is a Jordan triple
derivation, we have the
Peirce rules
$$ \{V_i, V_j, V_k\} \subeq V_{i-j+k}, \leqno(1.1) $$
which imply in particular that each space $V_j$ is a Jordan subtriple.
In addition, we have
$$ V_0 \square V_2 = V_2 \square V_0 = \{0\}. \leqno(1.2)  $$

The Jordan triple $V$ also carries a Jordan algebra structure, denoted
$V^{(e)}$, given by
$$ ab := L(a).b := \{ a,e,b\}. $$
Then $e$ is an idempotent in $V^{(e)}$ because $ee = \{e,e,e\} = e$, and
the Peirce decomposition of $V$ with respect to the tripotent $e$
coincides with
the Peirce decomposition of the Jordan algebra $V^{(e)}$ with respect
to the idempotent~$e$.

The multiplication operators in $V^{(e)}$ are given by $L(a) = a \square e$, so
that $L(e)\res_{V_2} = \id_{V_2}$ implies that $(V_2, e)$ is a unital Jordan
subalgebra of $V^{(e)}$. For the quadratic representation in $V^{(e)}$ we have
$$ P(e) = 2 L(e)^2 - L(e^2) = 2 L(e)^2 - L(e) = (2 L(e) - \1) L(e), $$
so that $P(e) = Q(e)^2$ vanishes on $V_0 \oplus V_1$ and restricts to the
identity on $V_2$. It follows in particular
that $(V_2, e, Q(e))$ is an involutive Jordan algebra
(cf.\ [Lo77, Th.~3.13]).

\msk\nin {\bf I.3. Orbits in $\oline{\cal D}$.}
Two tripotents $e,f \in V$ are said to be {\it orthogonal} if
$f \in V_0(e)$. In view of the Peirce rules (1.2), this
implies $\{f,f,e\} = \{e,f,f\} = (e \square f).f = 0$, so that we also have
$e \in V_0(f)$, i.e., orthogonality is a symmetric relation.
If this is the case, then $e + f$ also is a tripotent because the relations
$e \square f = f \square e= 0$ lead to
$$ \{e + f, e + f, e+f\} = \{e,e,e+f\} + \{f,f,e+f\}
= \{e,e,e\} + \{f,f,f\} = e + f. $$
We call a non-zero tripotent $e$ {\it primitive} if it cannot be
written as a sum of
two non-zero orthogonal tripotents and $e$ is said to be {\it
maximal} if there is
no non-zero tripotent orthogonal to~$e$.
A maximal tuple $(c_1,\ldots, c_r)$ of mutually orthogonal primitive tripotents
is called a {\it Jordan frame in $V$} and $r = \rank {\cal D}$ is
called the {\it rank
of ${\cal D}$}.
We fix a Jordan frame $(c_1,\ldots, c_r)$.
For $k = 0,1,\ldots, r$ we  then obtain tripotents
$$ e_k := c_1 + \ldots + c_k, $$
where it is understood that $e_0 = 0$.

We recall that each bounded symmetric domain ${\cal D}$ can be decomposed in a
unique fashion as a direct product of indecomposable, also called
{\it irreducible},
bounded symmetric domains:
$$ {\cal D} = {\cal D}_1 \times \ldots \times {\cal D}_m. \leqno(1.3) $$
Then the connected group $G := \Aut({\cal D})_0$ satisfies
$$ G \cong G_1 \times \ldots \times G_m, \quad \hbox{ where } \quad
G_j := \Aut({\cal D}_j)_0. \leqno(1.4) $$
If ${\cal D}$ is irreducible, then $G$ has exactly
$r + 1$ orbits in the closure $\oline{\cal D}$ of ${\cal D}$ in $V$
and $e_0, \ldots, e_r$ form a set of representatives
(cf.\ [Sa80, Th.~III.8.7]). For $k = 0$ we have
$G.e_0 = {\cal D}$ and for $k = r$ we obtain the {\it Shilov boundary}
$G.e_r = S$ ([Sa80, Th.~III.8.14]). We define the {\it rank of $x \in
\oline{\cal D}$} by
$$ \rank x = k \quad \hbox{ for } \quad x \in G.e_k $$
and observe that the rank function is $G$-invariant and classifies the
$G$-orbits in $\oline{\cal D}$.

If ${\cal D}$ is not irreducible, then (1.3/4) imply that the orbit of
$x = (x_1,\ldots, x_m) \in \oline{\cal D} = \prod_{j=1}^m
\oline{{\cal D}_j}$ is
determined by the $m$-tuple
$$ (\rank x_1, \ldots, \rank x_m) \in \N_0^m. $$
Here $(0,\ldots, 0)$ corresponds to elements in ${\cal D}$ and
$(\rk {\cal D}_1, \ldots, \rk {\cal D}_m)$ to elements in the product set
$S = S_1 \times \ldots \times S_m$.

\msk

\nin {\bf I.4. Spectral decomposition and spectral norm.}
Let $K$ be the stabilizer of $0 \in {\cal D}$ in~$G$. Then $K$ acts
as a group of automorphisms on the Jordan triple $V$ and each
element $z \in V$ is conjugate under $K$ to an element in
$\span_\R \{ c_1, \ldots, c_r\}$. For
$k.z = \sum_{j = 1}^r \lambda_j c_j$ the number
$$ |z| := \max\{|\lambda_1|, \ldots, |\lambda_r|\} $$
is called the {\it spectral norm of $z$}. Then the elements
$\tilde c_j := k^{-1}.c_j$ are orthogonal tripotents with
$$ z = \sum_{j = 1} \lambda_j \tilde c_j, $$
which is the spectral decomposition of $z$.
The spectral norm $|\cdot|$ is indeed a norm
on $V$ with
$$ {\cal D} = \{ z \in V \: |z| < 1\}. \leqno(1.5) $$

The following theorem relates the holomorphic arc-components in
$\partial {\cal D}$ to the
tripotents in~$V$.

\Theorem I.5. {\rm([Lo77, Th.~6.3])} For each holomorphic
arc-component $A$ of $\partial {\cal D}$ there exists a tripotent
$e$ in $A$ such that
$$ A = A_e := e + {\cal D}_e, \quad \hbox{ where } \quad
{\cal D}_e := {\cal D} \cap V_0(e) $$
is a bounded symmetric domain in $V_0(e)$. The map
$e \mapsto A_e$
yields a bijection from the set of non-zero tripotents of $V$ onto the set
of holomorphic arc-components of $\partial{\cal D}$. The Shilov boundary
$S$ coincides with the set of maximal tripotents.

An element $x \in \oline{\cal D}$ is contained in $A_e$ if and only if
$$ e = \lim_{n \to \infty} Q(x)^n.x. \leqno(1.6) $$
\qed

\msk
\nin{\bf I.6. Conformal completion of $V$.}
Let $G_\C$ denote the universal complexification
of the connected real Lie group $G$ and
$\tau$ the anti-holomorphic involution of $G_\C$ for which $G$ is the
identity component of the fixed point group $G_\C^\tau$.
Then the Lie algebra $\g_\C$ of $G_\C$ has a faithful realization by polynomial
vector fields of degree $\leq 2$ on $V$, which leads to a $3$-grading
$$ \g_\C = \g_+ \oplus \g_0 \oplus \g_-, $$
where $V \cong \g_+$ is the space of constant vector fields,
$\g_0$ consists of linear vector fields, and
$\g_-$ is the set of quadratic vector fields corresponding to the maps
$z \mapsto Q(z).v = \{z,v,z\}$ for $v \in V$.
By construction of the triple product, the
vector fields $\xi_v$ correspond to elements of the real Lie algebra
$\g = \L(G)$,
which implies that $\tau$ maps the constant vector field $v$
to the quadratic vector field $z \mapsto -\{z,v,z\}$. Hence
$\tau$ reverses the grading of $\g_\C$, i.e., $\tau(\g_j) = \g_{-j}$
for $j \in \{+,-,0\}$.
The Jordan triple structure on $V \cong \g_+$ then satisfies
$$ \{ x,y,z\} = {1\over 2} [[x,\tau.y], z]. \leqno(1.7)  $$
The subgroups
$$ G^\pm := \exp \g_\pm \quad \hbox{ and } \quad
G^0 := \{ g \in G_\C \: (\forall j)\ \Ad(g)\g_j = \g_j\} $$
satisfy
$$ G^\pm \cap G^0 = \{\1\} \quad \hbox{ and } \quad
(G^\pm \rtimes G^0) \cap G^\mp = \{\1\}. $$
Therefore $P^\pm := G^\pm G^0 \cong G^\pm \rtimes G^0$ are subgroups
of $G_\C$, and we obtain an embedding
$$ V \into X := G_\C/P^-, \quad v \mapsto \exp v \cdot P^-,  $$
called the {\it conformal completion of $V$}.
The elements of $G^+$ act on $V \subeq X$ by translations
$$ t_v \: x \mapsto x + v \leqno(1.8) $$
because $\exp v \exp x P^- = \exp(v + x) P^-.$
We further have $\tau(G^\pm) = G^\mp$ and $\tau(G^0) = G^0$.

For $w \in V$ we write $\tilde t_w$ for the map $X \to X$ induced by
the element
$\exp(-\tau(w)) = (\tau(\exp w))^{-1}$. For
$v \in V$ the condition $\tilde t_w.v \in V$, where $V$ is considered
as a subset of $X$,
is then equivalent to the invertibility of
$$ \1 + \ad v \ad (-\tau.w) + {1\over 4} (\ad v)^2 (\ad \tau.w)^2
= \1 - \ad v \ad(\tau.w) + {1\over 4} (\ad v)^2 \circ \tau \circ (\ad \tau)^2
\circ \tau \leqno(1.9)$$
([BN05, Cor.~1.10]).
In view of (1.7), this is precisely the Bergman operator
$$ B(v,w) = \1 - 2 v \square w + Q(v) Q(w). $$
We further have in $V$ the relation
$$ \tilde t_w.v = B(v,w)^{-1}.(v - Q(v).w). \leqno(1.10) $$

\msk
\nin{\bf I.7. Quasi-invertibility and transversality.} A
pair $(x,y) \in  V$ is called {\it quasi-invertible}
if $B(x,y) \in \End(V)$ is invertible. We write $x \top y$ if $(x,y)$
is quasi-invertible
and say that $x$ is {\it transversal} to $y$. We write
$x^\top := \{ y \in V \: x \top y\}$ for the set of all elements in
$V$ transversal to~$x$.

In the Jordan algebra $V^{(y)}$ with the product
$ab := \{a,y,b\}$ we have $L(a) = a \square y$ and
$P(a) = Q(a) Q(y)$ ([N\O{}04, App.~A]), so that
$$ B(x,y) = \id_V - 2 L(x) + P(x), $$
and in the unital Jordan algebra
$V^{(y)} \times \R$ with the identity element $\1 := (0,1)$  we have
$$ \1 - 2 L(x) + P(x)
= P(\1,\1) - 2 P(\1,x) + P(x,x) = P(\1 - x), $$
i.e., the quasi-invertibility of $(x,y)$ is equivalent to the
quasi-invertibility of $x$ in the Jordan algebra $V^{(y)}$.

\msk
\nin{\bf I.8. The $\sL_2$-triple associated to a tripotent.}
Let $e \in V$ be a tripotent, $f := \tau(e)$, $h := [e, f]$ and
$\g_e := \span_\R \{ h,e,f\}$. Then
$$ [h,e] = 2 \{e,e,e\} = 2e
\quad \hbox{ and } \quad
[h,f] = \tau [\tau h, e] = - \tau [h,e] = - 2 \tau e = - 2f, $$
so that $\g_e \cong \sL_2(\R)$ is a $3$-dimensional subalgebra of $\g$ with
$\g_e^\tau = \R(e + f)$.

(a) The operator $\ad_V h = 2 e \square e$
is diagonalizable with possible eigenvalues $0,1,2$. The corresponding
eigenspace decomposition
$V = V_0 \oplus V_1 \oplus V_2$
is the Peirce decomposition of the Jordan algebra
$V^{(e)}$ with multiplication $ab := \{ a,e,b\}$ with respect to the
idempotent $e$, i.e., $2 L(e).v_j = j v_j$ for $j = 0,1,2$.


(b) We observe that $P(e) = 2 L(e)^2 - L(e^2) = (2 L(e) - \1) L(e)$.
For $\lambda \in \R$ we therefore have for
$$ \eqalign{
&\ \ \ \ B(e, (1- \lambda) e)
= B((1- \lambda) e,e)
= \1 - (1- \lambda) 2 e \square e + (1-\lambda)^2 Q(e)^2 \cr
&= \1 - (1- \lambda) 2 L(e) + (1-\lambda)^2 P(e)
= \1 - (1- \lambda) 2 L(e) + (1-\lambda)^2 (2 L(e) - \1) L(e) \cr} $$
the relation
$$ B(e, (1-\lambda)e)v_j = \lambda^j v_j, \quad j = 0,1,2. $$

(c) From $Q(e) = Q(Q(e)e) = Q(e)^3$ we conclude that the antilinear map
$Q(e)$ is diagonalizable over $\R$ with
eigenvalues in $\{1,0,-1\}$, so that $Q(e)^2 = P(e) = (2 L(e) - \1)
L(e)$ implies
that
$$\ker Q(e) = \ker P(e) = V_0 \oplus V_1. \leqno(1.11) $$
  From $V_0 \square V_2 = V_2 \square V_0 = \{0\}$ we obtain for $x,y \in V_0$:
$$ \eqalign{
&\ \ \ \ B(e+ x, e+ y).e = e - 2 (e + x) \square (e + y).e + Q(e+x)Q(e+y) e \cr
&= e - 2 e - 2 x \square y.e + Q(e+x)(Q(e).e + Q(y).e + 2\{e,e,y\})\cr
&= -e -2 (e \square y).x  + Q(e+x).e = -e + (Q(e).e + Q(x).e +
\{e,e,x\})=0. \cr} $$

\Theorem I.9. {\rm([Lo77, Th.~8.11])}
Let $e \in V$ be a tripotent and
$V^{(e)}$ the corresponding Jordan algebra with product
$ab = \{a,e,b\}$. Identifying $e \in V$ with an element of $\g_+$, the
partial Cayley transform corresponding to $e$ is defined by
$C_e := \exp\Big({\pi \over 4}(e - \tau.e)\Big) \in G_\C$, and in Jordan
theoretic terms it is given as a partially defined map on $V$ by
$$ C_e
= t_e \cdot B(e,(1- \sqrt 2)e) \cdot \tilde t_e. $$
In particular
$$ C_e^{-1}(V) \cap V
= \{ v \in V \: B(e,v) \in \GL(V)\} = e^\top.
\qeddis

In [Lo77] Loos writes $B(e,-e)^{\half}$ instead of $B(e, (1- \sqrt
2)e)$, which makes sense because
$$B(e, (1- \sqrt 2)e)^2 = B(e, (1 - 2)e) = B(e,-e)$$ is
diagonalizable and the eigenvalues
$1,\sqrt 2$ and $2$ of $B(e,(1-\sqrt 2)e)$ are positive (I.8).

\msk\nin{\bf I.10.} The preceding theorem implies in particular that
the condition for an element $x \in V$ to lie in the
domain of the Cayley transform is precisely the transversality
condition $e \top x$.
If $x_2$ is the Peirce component of $x$ in $V_2$, then
[Lo77, Prop.~10.3] says that $e \top x$ is equivalent to the invertibility of
$e - x_2$ in the unital Jordan algebra $(V_2, e)$.

\Definition I.11. A hermitian scalar product
$\la \cdot, \cdot \ra$ on
$V$ is said to be {\it associative} if for $x,y,z ,w\in V$ we have
$$ \la \{x,y,z\},w \ra = \la x, \{y,z,w\}\ra, $$
which is equivalent to
$$ (z \square y)^* = y \square z \quad \hbox{ for } \quad y,z \in V. $$

According to [Lo77, Cor.~3.16], a scalar product with this property is given by
$$ \la x, y \ra := \tr(x \square y), $$
and for $0 \not= x \in V$ the operator $x \square x$ is non-zero and
positive semidefinite. In this sense $(V,\{\cdot,\cdot,\cdot\})$ is a
{\it positive
hermitian Jordan triple}.
\qed

\Lemma I.12. Let $e \in V$ be a tripotent, $V_j := V_j(e)$ its Peirce spaces,
and $z \in V_0$ with $|z| \leq 1$.
Further let $f := \lim_{n \to \infty} Q(z)^n.z$ denote the unique tripotent
contained in the holomorphic arc-component of $z$.
Then $\phi(z) := Q(z + e)\res_{V_1} \: V_1 \to V_1$ is an antilinear operator
which is symmetric with respect to the real scalar product
$(z,w) := \Re \tr(z \square w)$,
   and for $z \in V_1$ we have $\phi(z)v = 2\{z,v,e\}$.

If $|z| < 1$, then $\phi(z) + \1$ is injective ($\1$ stands for
$\id_{V_1}$), and for
$|z| = 1$ its kernel is
$$ \Fix(-Q(e + f)) \cap V_1(f) \cap V_1(e). $$

\Proof. For $v \in V_1$ we have
$$ \phi(z)v = \{ z + e, v, z + e \} = Q(z)v + Q(e)v + 2 Q(z,e)v, $$
and $Q(e)v \in V_{4-1} = V_3 = \{0\}$ as well as
$Q(z)v \in V_{0-1} = V_{-1} = \{0\}$ by the Peirce relations (1.1), so that
$\phi(z)v = 2 \{z,v,e\}.$

According to [Lo77, Lemma~6.7], the operator $\phi(z)$ on $V_1$ is
symmetric with respect to the real scalar product
$(\cdot,\cdot)$ on $V_1$, hence diagonalizable over
$\R$ with real eigenvalues.

Let $v \in V_1$ be an eigenvector for $\phi(z)$ corresponding to the eigenvalue
$\lambda \in \R$, i.e., $Q(z+e).v = \lambda v$.
Inductively we get
$$ Q(Q(z + e)^n.(z + e)).v = \lambda^{2n +1} \cdot v $$
for all $n \in \N_0$ from
$$ \eqalign{ Q(Q(z + e)^n.(z + e)).v
&= Q(Q(z + e)Q(z + e)^{n-1}.(z + e)).v \cr
&= Q(z + e) Q(Q(z + e)^{n-1}.(z + e))Q(z + e).v \cr
&= Q(z + e) Q(Q(z + e)^{n-1}.(z + e)).\lambda v
=  \lambda Q(z + e).(\lambda^{2n-1}.v)=  \lambda^{2n+1}v. \cr} $$

Since the inclusion $V_0 \into V$ is isometric with respect to the
spectral norm ([Lo77, Th.~3.17]), we have
$$ e + z \in e + \oline{{\cal D}_e} = \oline{A_e} \subeq \oline{\cal D}, $$
and the limit
$f  = \lim_{n \to \infty} Q(z)^n.z$
is a tripotent in $V_0(e)$ (Theorem~I.5).

As a consequence of the Peirce relations (1.2),
we obtain
$$ Q(e + z).(e + z) = Q(e)e + Q(z) z = e + Q(z) z, $$
and inductively
$$ Q(e + z)^n.(e + z) =  e + Q(z)^n.z \to e + f. $$
Therefore
$$ \lim_{n \to \infty} \lambda^{2n+1}v
=  \lim_{n \to \infty} Q(Q(z + e)^n.(z + e)).v
=  Q(e + f).v, $$
and the existence of the limit implies that $|\lambda| \leq 1$.
If $|\lambda| < 1$, then $Q(e + f).v = 0$, and
otherwise $Q(e + f).v = \lambda v$. It follows in particular
that each eigenvector for
$Q(e + z)$ on $V_1$ also is an eigenvector of $Q(e + f)$.

Suppose that $|\lambda| = 1$.
As a consequence of the Peirce rules, the sum $e + f$ is a
Jordan tripotent (I.3), and from $Q(e + f).v = \lambda v$ and
$\ker Q(e + f) = V_0(e + f) \oplus V_1(e + f)$ (I.8),
we derive
$v \in V_2(e + f)$, so that $(e + f) \square (e + f)
= e \square e + f \square f$ implies that $v \in V_1(f)$.

On the other hand $Q(e + f)$ is an antilinear involution of
$V_2(e + f) \supeq V_1(e) \cap V_1(f)$. We conclude
that
$$ \ker(\phi(z) + \1) = \ker(\phi(f) + \1)
= \Fix(-Q(e + f)) \cap V_1(f) \cap V_1(e).
\qeddis

To classify the $G$-orbits of transversal pairs in $\oline{\cal D}$,
we need a more explicit description of the image
$$ {\cal D}^C := C_e({\cal D}) $$
of ${\cal D}$ under the partial Cayley transform $C_e$
in terms of the Peirce decomposition of $V$.
To this end, we introduce the following notation:

\Definition I.13. Let $e \in V$ be a tripotent.
\litem{(1)} $(V_2, e, Q(e))$ is a unital involutive Jordan algebra.
We write $v^* := Q(e)v$ for the involution on $V_2$ and observe that
$V_2 = E \oplus i E$ for $E := \{ v \in V \: v^* = v\}$.
In this sense
$$\Re v = \half(v + v^*) = \half(v + Q(e)v) $$
is the component of $v$ in the real form $E$ of $V_2$.
The real Jordan algebra $E$ is euclidean and
we write $E_+ := \{ a^2 \: a \in E\}$ for its closed positive cone.
For $v,w \in E$ we write $v > w$ for $v - w\in \Int(E_+)$ and
$v \geq w$ for $v -w \in E_+$.
\litem{(2)} For $z \in V_0$ we define the antilinear map
$$ \phi(z)\: V_1 \to V_1, \quad v \mapsto 2\{ e,v,z\} = Q(e + z).v $$
(Due to the different normalization, the factor $2$ not present in [Lo77]).
\litem{(3)} We also define a hermitian map
$$ F \: V_1 \times V_1 \to V_2, \quad (z,w) \mapsto \{z,w,e\} $$
with
$$ F(z,w)^* = F(w,z) \quad \hbox{ and } \quad
F(z,z) > 0 \quad \hbox{ for } \quad 0 \not= z \in V_1.$$
For $u \in V_0$ with $|u| < 1$ we further define a real bilinear map
$$ F_u(z,w) = F(z, (\1 + \phi(u))^{-1}.w), $$
where we recall from Lemma~I.12 that $\1 + \phi(u)$ is invertible.
\qed

In the following proposition the missing factor $\half$ in front of $F$,
compared to [Lo77],  is due
to our different normalization of the triple product.

\Proposition I.14. {\rm([Lo77, Th.~10.8])} We have
$$ {\cal D}^C =C_e({\cal D})
=  \{ v = v_2 + v_1 + v_0 \in V_2 \oplus V_1 \oplus V_0 \:
|v_0| < 1, \Re(v_2 - F_{v_0}(v_1, v_1)) > 0\}.
\qeddis

To determine the closure of ${\cal D}^C$,
we need the following lemma, because there might be elements
$x_0 \in \partial {\cal D} \cap V_0$ for which the operator
$\phi(x_0) + \1$ is not invertible.

\Lemma I.15. Let $F$ be a finite-dimensional euclidean vector space,
$(A_n)_{n \in \N}$ a sequence of positive definite operators on $F$
converging to $A$
and $(v_n)_{n \in \N}$ a sequence of elements of $F$ converging to $v$.
If the sequence $A_n^{-\half}v_n$ is bounded, then
$v \in \im(A).$

\Proof. Since $A$ is symmetric, we have $\im(A)= \ker(A)^\bot$.
Let $w \in \ker(A)$. We have to show that $\la v, w \ra = 0$.
Since the sequence $A_n^{-\half} v_n$ is bounded, it contains a convergent
subsequence, and we may thus assume that it converges to some $u \in F$.
Then we get
$$ \eqalign{ \la v, w \ra
&= \lim_{n \to \infty} \la v_n, w \ra
= \lim_{n \to \infty} \la A_n^{\half} A_n^{-\half} v_n, w \ra
= \lim_{n \to \infty} \la A_n^{-\half} v_n, A_n^{\half} w \ra
= \la u, A^{\half} w \ra = \la u, 0\ra = 0. \cr} $$
This completes the proof.
\qed

\Lemma I.16. For each element
$v = v_2 + v_1 + v_0 \in \oline{{\cal D}^C}$ we have
$v_1 \in \im(\1 + \phi(v_0)).$

\Proof. Let $(v^n)_{n \in \N}$ be a sequence in ${\cal D}^C$
converging to $v$ and
write $v^n_j$, $j = 0,1,2,$ for its Peirce components.

We pick a linear functional $f \in E^*$ in the interior of the dual cone
of $E_+$, so that $f(x) > 0$ holds for $0 \not= x \in E_+$, and
observe that this implies that
$$ (v,w) := f(\Re F(v,w)) $$
defines a real scalar product on $V_1$. The argument in [Lo77, p.10.6]
shows that for each $z \in V_0$ the operator $\phi(z)$ is symmetric
with respect to this scalar product. According to Lemma~I.12, all its
eigenvalues $\lambda$ satisfy $|\lambda| \leq 1$ and even
$|\lambda| < 1$ for $|z| <1$, so that
$\1 + \phi(z)$ is a positive semidefinite symmetric operator which is
positive definite for $|z| < 1$.

  From $v^n \in {\cal D}^C$ we get
$$ |v_0^n| < 1 \quad \hbox{ and } \quad \Re F_{v_0^n}(v_1^n, v_1^n) \leq
\Re v_2^n, $$
which implies that
$$ \eqalign{ f(v_2^n)
&\geq f(\Re F_{v_0^n}(v_1^n, v_1^n))
= f(\Re F(v_1^n, (\1 + \phi(v_0^n))^{-1} v_1^n)) \cr
&= (v_1^n, (\1 + \phi(v_0^n))^{-1} v_1^n)
= ((\1 + \phi(v_0^n))^{-\half} v_1^n, (\1 + \phi(v_0^n))^{-\half}
v_1^n).\cr} $$
Therefore the sequence  $(\1 + \phi(v_0^n))^{-\half} v_1^n$ in $V_1$
is bounded, and Lemma~I.15 implies that
$$ v_1 = \lim_{n \to \infty} v_1^n \in \im(\1 + \phi(v_0)).
\qeddis

\msk\nin{\bf I.17.} From the preceding lemma one easily derives an explicit
description of the closure of ${\cal D}^C$ because the operator
$(\1 + \phi(v_0))^{-1}$ is well-defined on $\im(\1 + \phi(v_0))$.
This leads to
$$ \oline{{\cal D}^C} =  \Big\{ v \in V\:
|v_0| \leq 1, v_1 \in \im(\phi(v_0) + \1),
\Re\big(v_2 - F(v_1, (\1 + \phi(x_0))^{-1}v_1)\big) \geq 0\Big\}. $$
Since we do not need this description in the following, we
leave the details of its verification to the reader.

\Theorem I.18. {\rm(Orbits of transversal pairs)}
Let ${\cal D}$ be an irreducible bounded symmetric domain, not
necessarily of tube type.
If $(x,y) \in \oline{\cal D}$ is a transversal pair with
$\rk x = k$, then there exists a $g \in G$ with
$g.(x,y) = (e_k,z)$ with
$$ e_k =  c_1 + \ldots + c_k \quad \hbox{ and } \quad
z = -(c_{j+1} + \ldots + c_k) + \sum_{l = k+1}^r \lambda_l c_l, \quad
-1 \leq \lambda_{k+1} \leq \ldots \leq \lambda_r \leq 1. $$

\Proof. Since ${\cal D}$ is irreducible,
$G$ acts transitively on the set of elements of rank $k$, so that
we may w.l.o.g.\ assume that $x = e := e_k$. We then have to show
that each $G_e$-orbits in $e^\top \cap \oline{\cal D}$
contains an element of the form
$$ -(c_{j+1} + \ldots + c_k) + \sum_{l = k+1}^r \lambda_l c_l, \quad
-1 \leq \lambda_{k+1} \leq \ldots \leq \lambda_r \leq 1. $$
We recall the notation from Definition~I.13. For $y > 0$ in $E$ we
then find with (I.7)
$$ B(e-y,e) = \id_V - 2 L(e - y) + P(e-y) = P(e - (e- y)) = P(y).
\leqno(1.12) $$

Let $Q := G_{A_e}$ denote the stabilizer of the
holomorphic arc-component $A_e$ of $e$ in $\partial {\cal D}$ (which
is a maximal parabolic subgroup of $G$). Then the group
$Q^C := C_e \circ Q \circ C_e^{-1}$ acts naturally on ${\cal D}^C =
C_e({\cal D})$
and we also put
$$ Q^C_e := C_e \circ G_e \circ C_e^{-1} \subeq Q^C, $$
where  $G_e$ is the stabilizer of $e$ in $G$.

  From [Lo77, Lemma~10.7] we now obtain
$$ Q^C = \{ t_b \circ t_{v + F(v,v)} \exp(2e \square v) P(y)
\exp(\xi_w) \cdot k \:
b \in iE, v \in V_1, 0 < y \in E, w \in V_0, k \in K_e\}, $$
where $K_e := \{ g \in G : g.0 = 0, g.e = e\} \subeq \Aut(V)_e$ is
the set of all automorphisms of the Jordan triple $V$ fixing $e$
and $P(y)$ is the quadratic representation of
the Jordan algebra $V^{(e)}$ (cf.~I.7).
{}From the proof of [Lo77, Thm.~9.15] and the description of the Lie algebra
$\L(Q^C)$ in [Lo77, Prop.~10.6] it follows that for
$b \in iE, v \in V_1, 0 < y \in E$ and $k \in K_e$ we have
$$ t_b \circ t_{v + F(v,v)} \exp(2e \square v) P(y) k \in Q^C_e. $$
Moreover, the explicit calculations in the proof of [Lo77, Th.~10.8] further
imply that the map
$$ V_0 \to A_e = e + ({\cal D} \cap V_0), \quad w \mapsto \exp(\xi_w).e  $$
is bijective and that the Cayley transform fixes each $\xi_w$. This
implies that
$$ Q^C_e = \{ t_b \circ t_{v + F(v,v)} \exp(2e \square v) P(y) \cdot k \:
b \in iE, v \in V_1, 0 < y \in E, k \in K_e\}. $$

We observe that for $v \in V_1$ the Peirce rules imply that
$e \square v$ is a nilpotent operator on $V$ mapping $V_j \to V_{j+1}$.
For $x = x_2 + x_1 + x_0 \in \oline{{\cal D}^C}$ the $V_1$-component of
$$ t_{v + F(v,v)} \exp(2e \square v).x $$
is given by
$$ x_1 + v + \phi(x_0).v, $$
and since $-x_1 \in \im(\1 + \phi(x_0))$ by Lemma~I.16,
there is a unique $v \in \im(\1 + \phi(x_0))$ with
$$ t_{v + F(v,v)} \exp(2e \square v).x \in V_2 \oplus V_0. $$
  From that we conclude that each $Q^C_e$-orbit in $V$ through an element
$y = y_2 + y_1 + y_0 \in \oline{{\cal D}^C}$
contains an element of the form
$$ x_2 + x_0 \quad \hbox{ with} \quad |x_0| \leq 1
\quad \hbox{ and } \quad \Re x_2 \geq 0. $$
Applying elements of the form $t_v$, $v \in i E$, we may further assume that
$x_2 \in E$, so that we have an element in $E_+ \times {\cal D}_e$.
  From the explicit description of $Q^C_e$ we derive that
the intersection of the orbit of $x_2 + x_0 \in E + V_0$
with $E + V_0$ contains the orbit of
$x_2 + x_0$ under the group $Q'' := P(E_+) K_e$.

The orbits of $Q''$ on the set $E_+ \times \oline {\cal D}_e$ are
products of orbits
of the automorphism group
$G(E_+)$ of the symmetric cone $E_+$ in  $E$ and
orbits of the identity component of the group $K_e$ on
${\cal D}_e$. Since the action of the group $K_e$ preserves
the Peirce decomposition, it acts on ${\cal D}_e \subeq V_0$ as a subgroup of
$\Aut(V_0)$. The identity component of the latter group is obtained by
exponentiating elements of the Lie subalgebra
$V_0 + \tau(V_0) + [V_0, \tau(V_0)] \subeq \g_\C$
(here we use that ${\cal D}_e
= {\cal D} \cap V_0$ is an irreducible bounded symmetric domain; cf.~Th.~I.5),
and all the elements of this subalgebra commute with
the element $e \in V_2$ by the Peirce rules (I.2). Hence the image of
$K_e$ in $\Aut(V_0)$ contains the identity component of $\Aut(V_0)$.

For $e = e_k = c_1 + \ldots + c_k$, the orbits of $G(E_+)_0$, which
coincide with the orbits of the
full group $G(E_+)$, are represented by the elements
$$e_0 = 0, e_1 = c_1, \ldots, e_j = c_1 + \ldots + c_j, \ldots, e_k = e $$
([FK94, Prop.~IV.3.2]).
Since $(c_{k+1}, \ldots, c_r)$ is a Jordan frame in $V_0$, each orbit
of $\Aut(V_0)_0$ in
$V_0$ contains an element of the form
$$ \sum_{l = k+1}^r \lambda_l c_l, \quad
\lambda_{k+1} \leq \ldots \leq \lambda_r $$
(cf.\ [FK94, Prop.~X.3.2]).

Next we transfer this information back to the bounded picture, i.e.,
to $G_e$-orbits in $\oline{\cal D}$.
According to [Lo77, Prop.~10.3], we have
$$ C_e(x_2 + x_0) = C_e(x_2) + x_0 = (e + x_2)(e- x_2)^{-1} + x_0
\quad \hbox{ for } \quad x_2 \in V_2, x_0 \in V_0. \leqno(1.13) $$
For $e_j = c_1 + \ldots + c_j$, $j \leq k$, the element $e + e_j$ is
invertible in $V_2$, and we obtain
for $\tilde e_j := (e_j - e)(e_j + e)^{-1} = - C_e(-e_j) = C_e^{-1}(e_j)$ that
$C_e(\tilde e_j) = e_j$. An explicit calculation in the
associative Jordan algebra generated by
$c_1, \ldots, c_k$ quickly shows that
$$ \tilde e_j = - (e - e_j) = - e + e_j = - c_{j+1} - \ldots - c_k. $$
This completes the proof.
\qed

\msk\nin{\bf I.19.} For
the special case $k =  r$, i.e., $e \in S$, we have $V_0 = \{0\}$, so that
${\cal D}^C$ is the Siegel domain
$$ {\cal D}^C = \{ v = v_2 + v_1  \in V_2 \oplus V_1 = V\:
\Re(v_2 - F(v_1, v_1)) > 0\} $$
of type II.
In this case the orbits of $Q''_e$ are represented by elements of the form
$- e + e_j$, $j = 0,\ldots, r$, so that we obtain only finitely many orbits.
Observe that $\rk(-e + e_j) = r - j$, so that, even if $Q''$ is not connected,
it cannot have less orbits in $e^\top$ than its identity component.

\msk
There would be no substantial gain in the proof of Theorem I.18 by
assuming that $\cal
D$ is of tube type. However, in the sequel we will need only a special case of
the theorem, for which an easy direct proof (independant of the proof
of Theorem I.18)
can be offered.

\Lemma I.20. Suppose that ${\cal D}$ is irreducible and of tube type,
let $x\in S$ and $z\in \overline{\cal D}$, and assume
that
$x\top z$.  There exists $g \in G$ and
an integer $k, 0\leq k\leq r${\footnote
{$^{(\dagger)}$}{If $k=r$, use the convention that $ \sum_{j = r+1}^r
c_j= 0$.}}
such that
   $$ g(x) = e_r \quad \hbox{ and } \quad
g(z) = -\sum_{j=k+1}^rc_j = e_k - e_r\ .$$

\Proof. As $G$ is transitive on $S$, there is no restriction in assuming that
$x = e := e_r$. Now the transversality condition is equivalent to
$z$ belonging to the domain $V^\times + e$ of the Cayley transform
$C(z) := C_e(z) := (e + z)(e- z)^{-1}$ (cf.~(1.13)). Set $\zeta
= C(z)$ (Theorem~I.9). Then
$\zeta\in E_+ + i E$. The point $e$ is sent by the Cayley transform
``to infinity'', in such a way that the stabilizer of $e$ in $G$
corresponds via
conjugation by the Cayley transform to  a subgroup of the affine group of
$E^\C$, denoted by
$Q^C_e$, namely the semi-direct product of the translations by an
element of $iE$
and the group
$G(E_+)$ (after complexification to $E^\C$ of its action on $E$). By
using a translation, we see that in the $Q^C_e$-orbit of $\zeta$, there is an
element  of the form $\eta\in E_+$. Since ${\cal D}$ is irreducible, the
$G(E_+)$-orbits in
$E_+$ are known to be exactly the $r+1$ orbits of the elements
$e_k = \sum_{j=1}^kc_j$, with $k=0,1,\dots r$ (see [FK94,
Prop.~IV.3.2]). But now the
inverse Cayley transform of the element $\displaystyle \sum_{j=1}^kc_j $ is the
element
$\displaystyle e_k - e = -\sum_{j=k+1}^r c_j$. Hence the result.
\qed

\sectionheadline{II. Transversality and faces}

In this section we keep the notation from Section~I. In particular
${\cal D}$ is a circular irreducible bounded symmetric domain of
rank $r$ in $V$. The main result of this section
is that transversality of two elements
$x,y \in \oline{\cal D}$ is equivalent to the geometric property that
$x$ and $y$ do not lie in a proper face of the compact convex set
$\oline{\cal D}$ (Theorem~II.12).

\Definition II.1. (a) We call a non-empty convex subset $F$ of a
convex set $C$ a {\it face}
if for $0 < t < 1$ and $c,d \in C$ the relation
$t c + (1 - t) d \in F$ implies $c,d \in F$.
We write ${\cal F}(C)$
for the set of non-empty faces of $C$.
A face $F$ is called {\it exposed} if there exists a linear functional
$f \: V \to \R$ with
$$ F = f^{-1}(\max f(C)). $$
An {\it extreme point} $e \in C$ is a point for which $\{e\}$ is a face, i.e.,
$t c + (1 - t) d  = e$ for $c,d \in C$ and $0 < t < 1$ implies $c = d = e$.
We write $\Ext(C)$ for the set of extreme points of $C$.

The set of all faces of $C$
has a natural order structure given by set inclusion
whose maximal element is $C$ itself. All extreme points of $C$ are minimal
elements of this set, but $C$ need not have any extreme points.

Obviously, the intersection of any family of faces is a face.
We thus define for a subset $M \subeq C$ the
{\it face generated by $M$} by
$$ \Face(M) := \bigcap \{ F \subeq C \: F \in {\cal F}(C), M \subeq F\}. $$

(b) For a convex set $C$ in the vector space $V$ we write
$$ \algint(C) := \{ x \in C \: (\forall v \in C-C) (\exists \eps> 0) \
x + [0,\eps] v \subeq C\} $$
for its {\it algebraic interior}. If $V$ is finite-dimensional, then
$\algint(C)$ is the interior of $C$ in the affine subspace it generates.
\qed

\Remark II.2. (a) Suppose that $C$ is a convex subset of a finite-dimensional
vector space having non-empty interior. Then all proper faces of $C$
are contained in the boundary $\partial C$ and, conversely, the Hahn--Banach
Separation Theorem implies that each boundary point is contained in a
proper exposed face.

(b) For any non-empty convex subset of a finite-dimensional real vector space
the algebraic interior is non-empty. Hence, if $x$ belongs to the algebraic
interior of a face $F$, then $F$ is generated by  $\{x\}$.

(c) Since every face $E$ of a face $F$ of $C$ is also a face of $C$,
faces of exposed faces of $C$ are faces of $C$.
On the other hand,  every proper face is contained in an exposed face
(see (a)),
so that we obtain inductively, that for each face $F$ there exists a sequence
of faces
$$ F_0 = F \subeq F_1 \subeq \ldots \subeq F_n = C $$
for which $F_i$ is an exposed face of $F_{i+1}$ for $i = 0,\ldots, n-1$.
\qed

\Proposition II.3. The proper faces of the convex set
$\oline{\cal D}$ are the closures of the holomorphic arc-components in
$\partial {\cal D}$ and the Shilov boundary is the set of extreme points of
$\oline{\cal D}$.

In particular the group $G$ acts on the set
${\cal F}(\oline{\cal D})$ of faces of $\oline{\cal D}$.

\Proof.  For the fact that $S$ is the set of extreme points of
$\oline{\cal D}$ we refer to
[Lo77, Th.~6.5].

Next we use [Sa80, Lemma~III.8.11, Th.~III.8.13] to see that the proper
exposed faces $F$ of $\oline{\cal D}$ are the closures of the holomorphic
arc-components in $\partial {\cal D}$. Since the action of the group
$G$ on $\oline{\cal D}$ permutes the holomorphic arc-components in
$\partial {\cal D}$, it also permutes the exposed faces of~$\oline{\cal D}$.

We now claim that each face of $\oline{\cal D}$ is exposed.
Since every face $F$ of $\oline{\cal D}$ is generated by a suitable element
$x \in F$ (Remark~II.2), it suffices to show that the face generated by
any element $x \in \partial {\cal D}$ is exposed. Let
$A_x$ be the holomorphic arc-component of $\partial {\cal D}$
containing $x$. Then $\oline{A_x}$ is an exposed face of $\oline{\cal D}$
with $\algint(\oline{A_x}) = A_x$ (Theorem~I.5). Therefore the face
generated by $x$
coincides with $\oline{A_x}$, showing that every face of $\oline{\cal D}$
is exposed.
\qed

\Remark II.4. From the preceding proposition we know that the map
$F \mapsto \algint(F)$ is a $G$-equivariant bijection between the set
${\cal F}(\oline{\cal D})$ of faces of
$\oline{\cal D}$ and the set of holomorphic arc-components in $\oline{\cal D}$.

If ${\cal D}$ is irreducible, we
define the {\it rank of a face} by $\rk F := k$ if $\algint(F)$ consists of
elements of rank $k$. Since two holomorphic arc-components are conjugate under
$G$ if and only if their elements have the same rank (cf.\ Theorem~I.5),
the rank function
$$ \rk \: {\cal F}(\oline{\cal D}) \to \{0,\ldots, r\} $$
classifies the $G$-orbits in ${\cal F}(\oline{\cal D})$.
The stabilizer of a proper face, resp., a holomorphic arc-component
in $\partial {\cal D}$,
is a maximal parabolic subgroup of $G$ ([Sa80, Cor.~III.8.6]).

If ${\cal D} = {\cal D}_1 \times \ldots \times {\cal D}_m$ is a
direct product of the
irreducible domains ${\cal D}_j$, then each face $F$ of $\oline{\cal
D}$ is a product
$F_1 \times \ldots \times F_m$ of faces $F_j \in {\cal F}(\oline{\cal
D}_j)$, so that the
$G$-orbits in
$${\cal F}({\cal D}) \cong {\cal F}({\cal D}_1) \times \ldots \times
{\cal F}({\cal D}_m) $$
are classified by the $m$-tuple
$(\rk F_1, \ldots, \rk F_m)$.
\qed

In the following we shall prove that for
two elements $x,y \in \oline{\cal D}$ transversality
is equivalent to the {\it geometric transversality relation}
$\Face(x,y) = \oline{\cal D}$.
We start with the easy implication.

\Proposition II.5. If $x,y \in \oline{\cal D}$ are transversal, then
they are not contained
in a proper face, i.e., $\Face(x,y) = \oline{\cal D}.$

\Proof. If $x$ and $y$ are not geometrically transversal,
then $F := \Face(x,y)$ is a proper face of $\oline{\cal D}$,
hence of the form
$$ F = F_e = e + (\oline{\cal D} \cap V_0(e)) = (e + V_0(e)) \cap
\oline{\cal D} $$
for some tripotent $e\in V$ (Theorem~I.5, Prop.~II.3 and [Sa80,
Lemma~III.8.10] for the second equality).
Then $x,y \in F$ implies that $x,y \in e + V_0(e)$, so that I.8 leads to
$B(x,y).e = 0$. Thus $x$ and $y$ are not transveral. This proves the assertion.
\qed

\Example II.6. We consider the $r$-dimensional polydisc
$$ {\cal D} := \Delta^r := \{ z \in \C^r \: \max_j |z_j| < 1 \}
\subeq V = \C^r. $$
Let $(c_1, \ldots, c_r)$ denote the canonical basis of $\C^r$. The
corresponding Jordan triple structure is given by
$$ \{x,y,z\} = (x_1 \oline{y_1} z_1, \ldots, x_r \oline{y_r} z_r). $$
An element $z \in \C^r$ is a tripotent if $|z_j|^2z_j = z_j$ holds
for each $j$, which
means that either $z_j = 0$ or $|z_j| = 1$. We have
$$ \rk z = |\{ j \: |z_j| = 1\}|, $$
and the tripotents of maximal rank form the
$n$-dimensional torus $S = \T^n$, the Shilov boundary of $\Delta^r$.

Since the faces of $\oline{\cal D} = \oline\Delta^r$
are cartesian products of faces of the closed
unit disc
$$\oline\Delta = \{ z \in \C \: |z| \leq 1\}, $$
each face $F \in {\cal F}(\oline{\Delta^r})$
is a product $F_1 \times \cdots \times F_r$ of closed unit discs and
points in the
boundary of $\Delta$. For a subset $M \subeq \oline\Delta^r$, it
follows that the face
generated by $M$ is given by
$$ \Face(M) = F_1 \times \cdots \times F_r, \quad
F_j = \cases{
\{s\} & if $m_j = s \in \partial \Delta$ for all $m \in M$ \cr
\Delta & otherwise. \cr} $$
It follows in particular that $x,y \in \oline{\cal D}$ are contained
in a proper face
if and only if $x_j = y_j \in \partial \Delta$ holds for some $j$.

For $k \leq r$ we consider the tripotent
$e_k := c_1 + \ldots + c_k$. Then
$$ V_2 = \C^k \times \{0\}^{r-k} \quad \hbox{ and } \quad V_0 =
\{0\}^{k} \times \C^{r-k}. $$
An element $x \in \oline{\Delta^r}$ is transversal to $e_k$ if and only if
$e_k - (x_1,\ldots, x_k,0,\ldots, 0)$
is invertible in the unital Jordan algebra $(V_2, e_k)$, which means that
the first $k$ components of $x$ are different from~$1$ (I.10).
That this is not the case means that one component $x_j$, $j \leq k$,
equals $1$, and therefore $\Face(e_k,x) \not= \oline{\cal D}$.
If, conversely, $\Face(e_k,x) \not= \oline{\cal D}$, then $e_k,x$ are
contained in a proper
face of $\oline{\Delta^r}$ which implies that $x_j = 1$ for some $j \leq k$.
\qed

\Proposition II.7. Let $e \in V$ be a tripotent,
$V = \sum_{j = 0}^2 V_j$ the corresponding Peirce decomposition and
$p_j \: V \to V_j$ the projection along the other Peirce components.
Then each $V_j$ is a positive hermitian Jordan triple and we have
$$ {\cal D}_j = V_j \cap {\cal D} = p_j({\cal D}). $$
In particular, each map $p_j \: V \to V_j$ is a contraction with
respect to the spectral norms
determined by the domains ${\cal D}$ and ${\cal D}_j$.

\Proof. Let $\la \cdot, \cdot \ra$ be an associative hermitian scalar
product on
$V$ (Definition~I.11). Then the Peirce decomposition is orthogonal
with respect to
$\la \cdot, \cdot \ra$, so that it provides an orthogonal decomposition of
$V$ into $3$ Jordan subtriples ([Lo77, Th.~3.13]).

Clearly the restriction of the scalar product to each $V_j$ provides
an associative
scalar product on $V_j$ and for each $v \in V_j$ the operator $v
\square v$ is positive
semidefinite
on $V,$ which implies in particular that its restriction to $V_j$ is
positive semidefinite.
Hence each $V_j$ is a positive hermitian Jordan triple.

According to [Lo77, Th.~3.17], the inclusion maps $V_j \into V$ are
isometric with respect to the
spectral norm, which means that
$$ {\cal D}_j = V_j \cap {\cal D} = \{ z \in V_j \: |z| < 1\} $$
holds for the corresponding bounded symmetric domains.

To see that the projections $p_j$ are contractive with respect to the 
spectral norm,
let
$v \in V$ and $v_j = p_j(v)$ its component in $V_j$. For each unit vector
$w \in V_j$ the orthogonality of the Peirce decomposition implies that
$$ \la v \square v.w, w \ra
= \sum_{k,l= 0}^2 \la v_k \square v_l.w, w \ra
= \sum_{k= 0}^2 \la v_k \square v_k.w, w \ra
\geq \la v_j \square v_j.w, w \ra, $$
which leads for the spectral norm $|v_j|$ to
$$ \eqalign{ |v_j|^2
&= \|v_j \square v_j\|_{V_j}
= \sup \{ \la v_j \square v_j.w, w \ra \: w \in V_j, \la w,w \ra = 1\}\cr
&\leq \sup \{ \la v \square v.w, w \ra \: w \in V_j, \la w,w \ra = 1\}
\leq \sup \{ \la v \square v.w, w \ra \: w \in V, \la w,w \ra = 1\}
= |v|^2. \cr} $$
Since the inclusion $V_j \into V$ is isometric,
$p_j$ is a contraction with respect to the spectral norm, and therefore
${\cal D}_j \subeq p_j({\cal D}) \subeq {\cal D}_j$
proves equality.
\qed

\Corollary II.8. If $F$ is a proper face of $\oline{\cal D}_j$, then
$p_j^{-1}(F)$ is a proper face of $\oline{\cal D}$.
\qed

\Definition II.9. Suppose that $e \in V$ is a tripotent with $V_2(e)
= V$, so that $Q(e)$
is an antilinear involution on $V$ turning $(V, e, Q(e))$ into an
involutive unital Jordan algebra.
As in Section~I, we endow $V$ with the spectral norm $|z|$ whose open
unit ball is ${\cal D}$.

A {\it state of the unital involutive Jordan algebra $V$} is a linear
functional $f \: V \to \C$ with
$$  1 = f(e) = \|f\| := \sup |f({\cal D})|.
\qeddis

\Remark II.10. If $f$ is a state on $V$ and $y \in \oline {\cal D}$ with
$f(y) = 1$, then $e$ and $y$ lie in the proper face
$\{ z \in \oline{\cal D} \: \Re f(z) = 1\}$.
\qed

\Proposition II.11. If $y \in \oline{\cal D}$ and $e - y$ is not invertible
in the unital Jordan algebra $(V,e)$, there exists a state $f$ of $V$
with $f(y) = 1$.

\Proof. We endow $V$ with the associative scalar product
$\la z,w\ra := \tr(z \square w)$ (cf.~Def.~I.11).

By assumption $e - y$ is not invertible, which implies that the left
multiplication
$L(e - y) = (e - y) \square e$ is not invertible. Pick $v \in \ker L(e - y)$
with $\la v, v \ra = 1$. We consider the linear functional
$$ f \: V \to \C, \quad f(z) := \la L(z).v,v \ra $$
satisfying
$f(e) = \la v, v \ra = 1$ and
$$f(y) = \la L(y).v, v \ra = \la L(e).v, v \ra = f(e) = 1.$$

It remains to show that $f$ is a state.
Let $E := \{ z \in V \: z^* = Q(e)z = z\}$ denote the euclidean Jordan algebra
with $V \cong E \otimes_\R \C$ and unit element $e$.
We write $E_+$ for the closed positive cone in $E$. This is the set
of all those elements
$z$ for which there exists a system $c_1,\ldots, c_k$ of orthogonal
idempotents with
$e = c_1 + \ldots + c_k$ and non-negative real numbers $\lambda_j$ with
$$ z = \sum_{j = 1}^k \lambda_j c_j. $$
For such elements $z \in E_+$ we then have
$$ f(z)
= \sum_{j = 1}^k \lambda_j \la L(c_j).v, v \ra
= \sum_{j = 1}^k \lambda_j \la c_j \square c_j.v, v \ra \geq 0 $$
because $L(c_j) = c_j \square e = c_j \square c_j$ follows from
$c_j \square (e - c_j) = 0$ (I.2) and
the operators $c_j \square c_j$ are positive semidefinite on $V$
([Lo77, Cor.~3.16]).
We conclude that $f(E) \subeq \R$, so that
$f(z^*) = \oline{f(z)}$ for all $z \in V$.

  From $Q(e)^{-1} = Q(e)$ we derive
$Q(Q(e).z) = Q(e) Q(z) Q(e) = Q(e) Q(z) Q(e)^{-1}$, so that
$Q(e) \: z \mapsto z^*$
is a Jordan triple automorphism of $V$, hence an isometry for the spectral norm
$|\cdot|$ on $V$. This implies that $Q(e){\cal D} = {\cal D}$ and
therefore that for
$z = x + iy \in {\cal D}$, $x,y \in E$, we have
$$ |x| = \half |z + z^*| \leq \half(|z| + |z^*|) = |z|. $$
For the map $\Re \: V\to E, z \mapsto \half(z + z^*)$ this means that
${\cal D}_E := {\cal D} \cap E = \Re({\cal D}).$

For the functional $f$ we thus obtain
$$ \|f\|
= \sup |f({\cal D})|
= \sup \Re f({\cal D})
= \sup f(\Re {\cal D}) = \sup f({\cal D}_E). $$
In view of the Spectral Theorem for euclidean Jordan algebras ([FK94]), we have
$$ {\cal D}_E = (e - E_+) \cap (-e + E_+) \subeq e - E_+, $$
so that $f(z) \geq 0$ for $z \in E_+$ leads to
$\|f\| = \sup f({\cal D}_E) = f(e) = 1.$
This means that $f$ is a state.
\qed

\Theorem II.12. Two elements $x,y \in \oline{\cal D}$ are transversal
if and only if they are not contained in a proper face, i.e.,
$$ x \top y \quad \Longleftrightarrow \quad \Face(x,y) = \oline{\cal D}. $$

\Proof. In view of Proposition~II.3, geometric transversality is also invariant
under the action of the group $G$. On the other hand transversality
is invariant
under $G$ ([C\O{}01]), so that it suffices to assume that
$x = e$ is a Jordan tripotent. In view of Proposition~II.5, it
suffices to show that if
$e$ is not transversal to $y \in \oline{\cal D}$, then both $e$ and
$y$ lie in a proper face of
$\oline{\cal D}$.

For $e = 0$ we have
$\Face(x,e) = \oline{\cal D}$ because $e \in {\cal D} =
\algint(\oline{\cal D})$ and also
$e \top x$ for all $x \in \oline{\cal D}$ because $B(x,e)= \id_V$.

We may therefore assume that $e \not=0$. We have to show that if $e$
and $y$ are not transversal,
then they are contained in a proper face of $\oline{\cal D}$.
That $y$ is not transversal to $e$ is equivalent to the element $e - y_2$ being
not invertible in the unital Jordan algebra $V_2(e)$ (I.10).
In view of Proposition~II.11, combined with Remark~II.10,
$e$ and $y_2$ are contained in a proper face $F$ of
the convex set $\oline{\cal D}_2$. Hence
$e$ and $y$ are contained in the proper face $p_2^{-1}(F)$ of $\oline{\cal D}$
(Corollary~II.8).
\qed

\Example II.13. Let $p,q\in \N$, $r := \min(p,q)$, and $\|\cdot\|$
denote the euclidean norm on
$\C^p$, resp., $\C^q$. On the matrix space
$V := M_{p,q}(\C) \cong \Hom(\C^q,\C^p)$ we write
$|X|$ for the corresponding operator norm. Then
$$ {\cal D} := \{ X \in M_{p,q}(\C) \: |X| < 1 \} $$
is a bounded symmetric domain.
The pseudo-unitary group
$\UU_{p,q}(\C)$ acts transitively on ${\cal D}$ by
$$ \pmatrix{a & b \cr c & d \cr}.z  := (az + b)(cz + d)^{-1}, $$
the effectivity kernel of this action is $\T \1$, so that
$G = \Aut({\cal D})_0\cong \PU_{p,q}(\C)$. The $3$-grading of $\g_\C$
is induced
by the $3$-grading of $\gl_{p+q}(\C)$ given by
$$ \gl_{p+q}(\C)_+ = \pmatrix{ 0 & M_{p,q}(\C) \cr 0 & 0 \cr}, \quad
\gl_{p+q}(\C)_0 = \pmatrix{ \gl_p(\C) & 0 \cr 0 & \gl_q(\C) \cr} $$
and
$$ \gl_{p+q}(\C)_- = \pmatrix{ 0 & 0 \cr M_{q,p}(\C) & 0 \cr}. $$
We further have
$$ \uu_{p,q}(\C) = \Big\{\pmatrix{ a & b \cr b^* & d\cr} \: a^* = -a,
d^* = -d\Big\}. $$
The vector field associated to the one-parameter group given by
$\exp\Big(t\pmatrix{a & b \cr c & d \cr}\Big)$ is given by
$z \mapsto b - az - zd - zcz,$ so that the Jordan triple structure on
$V = M_{p,q}(\C)$ satisfies $Q(z)(w) = zw^*z$, which leads to
$$ \{a,b,c\} = \half(ab^*c + cb^*a). $$
In particular the Bergman operator satisfies
$$ B(v,w)z
= z - 2 v \square w.z + Q(v) Q(w)z
= z - (vw^*z + zw^*v) + v(wz^*w)^*v
= (\1 - vw^*)z(\1 - w^*v). $$
  From that it follows that $v \top w$ is equivalent to the invertibility
of $\1 - w^*v$ in the algebra $M_q(\C)$.

An element $e \in M_{p,q}(\C)$ is a tripotent if and only if
$ee^*e = e$, which implies that $ee^*$ and $e^*e$ are orthogonal
projections, and that
$e$ defines a partial isometry $\C^q \to \C^p$. If
$K := \ker(e)$ and $R := \im(e)$, then the face $F_e$ of $\oline{\cal
D}$ consists of
all matrices $z \in \oline{\cal D}$ with $z.v= e.v$ for $v \in \ker(e)^\bot$.
For $k = \rank(e)$ and an orthonormal basis $v_1, \ldots, v_k$ of
$\ker(e)^\bot$ and
$w_i := e.v_i$, we have
$$ F_e = \{ z \in \oline{\cal D} \: (\forall i)\ \la zv_i, w_i \ra = 1\}. $$
  From this description of the faces of $\oline{\cal D}$ it follows that
an element $z \in \oline{\cal D}$ is contained in a proper face if and only if
its restriction to some one-dimensional
subspace of $\C^q$ is isometric, i.e., if and only if $|z| = 1$.
Two elements $z,w$ generate a proper face if and only if there exists
a unit vector
$v \in \C^q$ for which $z.v = w.v$ is a unit vector in $\C^p$.

A Jordan frame is given by the matrices
$c_j := E_{jj}$, $j =1,\ldots, r$, with a single non-zero entry $1$
in position $(j,j)$.
The rank of ${\cal D}$ is $r$ and $e_r := c_1 + \ldots + c_r$ is a maximal
tripotent
with
$$ S = G.e_r = \cases{
\{ z \in M_{p,q}(\C) \: z^*z = \1\} & if $q \leq p$ \cr
\{ z \in M_{p,q}(\C) \: zz^* = \1\} & if $p \leq q$. \cr} $$
For $q \leq p$ this is the set of isometries $\C^q \into \C^p$
and for $p \leq q$ this is the set of all adjoints of isometries
$\C^p \to \C^q$.

Let $e_k := c_1 + \ldots + c_k$ be the canonical tripotent of rank $k$.
Writing an element $z \in M_{p,q}(\C)$ as a block matrix
$$z = \pmatrix{z_{11} & z_{12} \cr z_{21} & z_{22} \cr} \quad \hbox{
with } \quad
z_{11} \in M_k(\C), z_{12} \in M_{k,q-k}(\C), z_{21} \in M_{p-k, k}(\C), z_{22}
\in M_{p-k,q-k}(\C), $$
we have
$$ 2\{e,e,z\}
= ee^*z + ze^*e
= \pmatrix{ \1 & 0 \cr 0 & 0 \cr} \pmatrix{z_{11} & z_{12} \cr z_{21}
& z_{22} \cr}
+ \pmatrix{z_{11} & z_{12} \cr z_{21} & z_{22} \cr}\pmatrix{ \1 & 0
\cr 0 & 0 \cr}
= \pmatrix{2 z_{11} & z_{12} \cr z_{21} & 0\cr}. $$
This shows that
$$ V_2(e_k) \cong M_k(\C), \quad
V_1(e_k) \cong M_{k,q-k}(\C) \oplus M_{p-k,k}(\C) \quad \hbox{ and } \quad
V_0(e_k) \cong M_{p-k,q-k}(\C), $$
and therefore
$$ F_e = \Big\{\pmatrix{\1 & 0 \cr 0 & z} \: z \in M_{p-k,q-k}(\C),
|z| \leq 1\Big\}. $$
For $k = r$ we see in particular that $V_0(e_r) = 0$.
\qed

\sectionheadline{III. Orbits of triples in the Shilov boundary}

In this section we obtain the key result for our classification of
triples in $S$
in the tube type case.  We show that if
$(c_1,\ldots, c_r)$ is a Jordan frame in
$E$, then each
$G$-orbit in $S\times S\times S$ meets the Shilov boundary $T \cong
\T^r$ of the corresponding
polydisc. We further show that the polydiscs arising in this result can also be
characterized directly as the intersections of ${\cal D}$ with
$r$-dimensional subspaces of $V$, or, equivalently, as isometric images of
polydiscs under affine maps $\C^r \to V$, mapping $\Delta^r$
isometrically into ${\cal D}$.
In particular we show that any such affine map is linear.

\Theorem III.1. Suppose that ${\cal D} \subeq V$ is of tube type,
$(c_1, \ldots, c_r)$ is a Jordan frame in $V$, and
$$T := S \cap \span \{c_1, \ldots, c_r\}
= \Big\{ \sum_{j=1}^r \lambda_j c_j \: (\forall j)\ |\lambda_j| =1\Big\} $$
is the corresponding $r$-torus in $S$.
Then for each triple $(e,f,g) \in S$ there exists a $g \in G$ with
$g.e, g.f, g.h \in T$.

\Proof. Since Jordan frames and $G$ decompose according to the decomposition
of ${\cal D}$ into products of irreducible domains, it suffices to prove the
assertion for irreducible domains. We prove the assertion by induction on the
rank $r$ of ${\cal D}$.  Observe that the algebraic interior of any face
$F$ of
$\cal D$ is a bounded symmetric space of tube type. In fact, let $E$
be a euclidean
Jordan algebra which has $V$ as its complexification. Let $F$ be a face of
rank $k$. Then $F$ contains a tripotent $c$ of rank $k$ and there
exists a Jordan
frame $(c_1,c_2,\dots,c_r)$ in $E$ such that $ c=\sum_{j=1}^k\lambda_j c_j$,
with $ |\lambda_j|=1$ for $1\leq j\leq k$. Then
$V_0(c)$ is the complexification of the euclidean Jordan algebra
$E_0(c) = E_0(c_1+c_2+\dots+c_k)$. For  $z\in V_0(c)$, the spectral 
norm relative to
$V_0(c)$ coincides with the spectral norm in $V$, and so $V_0(c)\cap
{\cal D}={\cal
D}_0$ is the bounded symmetric domain of tube type associated to the euclidean
Jordan algebra $E_0(c)$. As
$$ \algint(F) = \algint(F_c) =c+({\cal D}\cap V_0(c))=c+{\cal D}_0, $$
we see that $\algint(F)$ is a bounded symmetric domain of tube type.

{\bf Case 1:} If $\Face(e,f,h)$ is proper, then its algebraic interior
is a bounded symmetric domain of tube type ${\cal D}'$ of smaller
rank and $(e,f,h)$ are contained in its Shilov boundary. In fact,
according to Theorem~I.5 and Proposition~II.3,
for each face $F$ of $\oline {\cal D}$ corresponding to the
holomorphic arc-component $A = \algint(F)$, the Shilov boundary
of $A$ is given by
$$ S_A = \Ext(\oline A) = \Ext(F) = \Ext(\oline{\cal D}) \cap F = S \cap F. $$
Since every element of $\Aut({\cal D'})_0$ is the restriction of an element
of $\Aut({\cal D})$ ([Sa80, \break Lemma~III.8.1]),
in this case the result follows from the induction hypothesis if $r > 1$.
If $r = 1$, then each proper face of $\oline{\cal D}$ is an extreme point,
so that the assumption that $e,f,h$ lie in a proper face implies
$e = f = h$. In this case we further have $c_1 \in S$, so that the
assertion follows
from the transitivity of the action of $G$ on $S$.

{\bf Case 2:} We assume that some pair $(e,f)$, $(f,h)$ or $(e,h)$ is
transversal.
We may w.l.o.g.\ assume that $(e,f)$ is transversal.
Then $\Face(e,f,h) \supeq F(e,f) = \oline{\cal D}$ by Theorem~II.12, and
$G.(e,f)$ contains $(e,-e)$ because $\rk f = \rk e = r$
(Lemma~I.20). Therefore the orbit of $(e,f,h)$ contains an
element of the form $(e,-e,h)$. Now the assertion follows from the
Spectral Theorem for unitary elements in $V$ (cf.\ [FK94, Prop.~X.2.3])
and (A.5) in the appendix.


{\bf Case 3:} $\Face(e,f,h) = \oline{\cal D}$, but
neither $(e,f)$, nor $(f,h)$ or $(e,h)$ is transversal.
Since $G$ acts transitively on $S$, we may w.l.o.g.\ assume that
$e = e_r = c_1 + \ldots + c_r$.
Consider the proper face $F := \Face(f,h)$ of $\oline{\cal D}$.
Then we have
$$ \oline{{\cal D}} = \Face(e,f,h) = \Face(\{e\} \cup F), $$
and for any $x \in \algint(F)$ we obtain
$$ \oline{{\cal D}} = \Face(\{e\} \cup F) = \Face(e,x), $$
which means that $e$ and $x$ are transversal (Theorem~II.12).

Now we need the classification of $G$-orbits in the set of
transversal pairs, which shows that the pair $(e,x)$ is conjugate to
an element of the
form $(e, -e + e_j)$ (Lemma~I.21). The face
$$ \eqalign{ F'
&= \Face(-e + e_j) = - \Face(e - e_j)
= -(e - e_j) + (V_0(e - e_j)\cap {\cal D})
= (e_j - e) + (V_2(e_j) \cap {\cal D}) \cr} $$
is a bounded symmetric domain of tube type of rank $j$, and
$(e,f,h)$ is conjugate to a triple of the form
$(e, f', h')$ where $f', h'$ are two elements in the Shilov boundary of $F'$,
where they are transversal because they generate $F'$ as a face
(Theorem~II.12).
Next we observe that the Peirce rules imply that by exponentiating
elements of the centralizer of $e- e_j$ in $\g$ we generate
the identity component $G^0$ of the
group $\Aut({\cal D} \cap V_0(e-e_j))$ and its elements $g$ act on
$e_j - e + z$ by
$$  g.(e_j - e + z) = (e_j - e) + g.z $$
because they commute with the translation $t_{e_j - e}$.
Now we conclude the proof by applying the special case of transversal
elements which has already been taken care of, to see that the $G^0$-orbit
of $(e,f',h')$ intersects $T$.
\qed

\Remark III.2. If ${\cal D}$ is not of tube type, then the Cayley
transform $C = C_e$ leads to a
realization of ${\cal D}$ as a Siegel domain ${\cal D}^C$ of type~$II$, and
since $C_e(-e) = 0$, the stabilizer $G_{e,-e}$ of $\pm e$ in $G$
corresponds to the stabilizer $Q^C_{e,-e} := C_e(G_{e,-e})$
of $0$ in the affine group $Q^C_e$,
and the identity component of this group is $G(E_+)_0 K_e$
(see the proof of Theorem~I.18).
The Shilov boundary of ${\cal D}^C$
is the set
$$\{ (v_2, v_1) \in V = V_2 \oplus V_1 \: \Re v_2 =  F(v_1, v_1)\}, $$
and from this description it is clear that no element $v_2 + v_1$
with $v_1 \not= 0$
is conjugate under $Q^C_{e,-e}$
to an element in $\span_\R \{ c_1,\ldots, c_r\} \subeq V_2$. Therefore
the condition that ${\cal D}$ is of tube type is necessary for the
conclusion of
Theorem~III.1.
\qed

\Example III.3. The simplest example of a bounded symmetric domain
not of tube type
is the matrix ball ${\cal D} \subeq \C^{n}$ for $n > 1$. Its rank is
$r = 1$ and in this case $G \cong \PSU_{n,1}(\C)$ (cf.\ Example~II.13).

To $z \in {\cal D}$ we assign the one-dimensional subspace
$L_z := \C\pmatrix{z \cr  1\cr} \in \C^{n+1}$. Endowing $\C^{n+1}$
with the indefinite hermitian
form $h$ given by
$$ h(z,w) := z_1 \oline{w_1} + \ldots + z_n \oline{w_n} - z_{n+1}
\oline{w_{n+1}}, $$
we see that ${\cal D}$ corresponds to the set of lines on which $h$ is
negative definite, and its Shilov boundary, the sphere $S \cong \SS^{2n-1}$,
corresponds to the set of isotropic lines. In this picture the action of
$\SU_{n,1}(\C)$ on ${\cal D}$ comes from the natural action of this
group on the
one-dimensional subspaces of $\C^{n+1}$.

Fixing a unit vector $e \in S$, the pair $(e,-e)$ corresponds to two
different isotropic lines $L_e$ and $L_{-e}$ in $\C^{n+1}$, and the stabilizer
of this pair in $\UU_{n,1}(\C)$ fixes the non-degenerate subspace
$L_e + L_{-e}$, and also its orthogonal complement of dimension
$n-1$. We conclude
that $\UU_{n,1}(\C)_{e,-e} \cong \R^\times \times \UU_{n-1}(\C)$,
and that no line $L_z \not\subeq L_e + L_{-e}$ can be moved by $\UU_{n,1}(\C)$
into the plane $L_e + L_{-e}$. On the other hand, the set of
isotropic lines in the plane
$L_e + L_{-e}$ corresponds to the circle in $S$ obtained by intersecting
$S$ with the boundary of a one-dimensional disc $\Delta \subeq {\cal
D}$ of size~$1$,
which in particular is a polydisc of maximal rank. This shows quite
directly that there are triples in $S$ that cannot be moved into the
one-dimensional
space $\C e$, so that Theorem~III.1 does not hold.

That Theorem~III.1 fails in this context, can be expressed
quantitatively by the observation
that
$$ F(\C v_1, \C v_2, \C v_3) :=
{h(v_1, v_2) h(v_2, v_3) h(v_3, v_1) \over
h(v_2, v_1) h(v_3, v_2) h(v_1, v_3)} $$
is a well-defined function on the set of triples of pairwise different
isotropic lines in $\C^{n+1}$ which is invariant under the
pseudo-unitary group $U_{n,1}(\C)$.
The function $F$ is  related to the {\it Cartan
invariant} (for a presentation and a generalization of this invariant we refer
to [Cl05]).
\qed

\Example III.4. The matrix ball
${\cal D} \subeq M_n(\C)$ is a symmetric domain of tube type with
Shilov boundary $S = U_n(\C)$, the unitary group.
The maximal polydiscs in ${\cal D}$ are obtained by intersecting ${\cal D}$
with the set of all matrices that are diagonal with respect to some
fixed orthonormal basis of $\C^n$ with respect to the standard scalar product.
A particular Jordan frame consists of the matrix units
$c_j := E_{jj}$, $j =1,\ldots, n$, whose span is the set of diagonal matrices.
Therefore Theorem~III.1 states that each triple
$(s_1, s_2, s_3)$ of unitary matrices can be diagonalized by
an element $g \in U_{n,n}(\C)$, acting on $U_n(\C)$ by
$$ \pmatrix{a & b \cr c & d \cr}.z = (a z + b)(cz+d)^{-1}. $$
The compact subgroup $U_n(\C) \times U_n(\C)$ acts linearly by
$(a,d).z = azd^{-1}$, and under this group each pair
$(s_1, s_2)$ is conjugate to a pair of the form
$(\1, s_2')$, where the stabilizer of $\1$ is the diagonal subgroup,
acting on the second component by $(a,a^{-1}).s_2 = as_2 a^{-1}$,
so that $s_2'$ can be diagonalized by conjugating with a suitable
element $a \in \UU_n(\C)$.
This means that diagonalizability of pairs reduces to classical linear algebra,
but diagonalizability of triples requires the non-linear action of
$U_{n,n}(\C)$
and Theorem~III.1.

A classification of the conjugation orbits of $U_n(\C)$ in
$U_n(\C)^2$ is given in [FMS04],
but since $U_n(\C)$ is much smaller than $U_{n,n}(\C)$, this
classification leads to
infinitely many orbits.
\qed

\subheadline{Polydisc in bounded symmetric domains}

Let ${\cal D} \subeq V$ be a bounded symmetric domain of rank $r$ and
$\Delta^r \subeq \C^r$ the $r$-dimensional unit polydisc.
We endow $\C^r$ with the metric defined by the sup-norm
$$ |z| := \max \{|z_1|, \ldots, |z_r|\} $$
and $V$ by the metric defined by the spectral norm, also denotes $|z|$.

\Theorem III.5. Any affine isometric map $f \: \C^r \to V$ mapping
$\oline\Delta^r$ into $\oline{\cal D}$ is linear
and preserves the rank, i.e., for each $x \in \oline\Delta^r$ we have
$$ \rk f(x) = \rk x. $$
Moreover, it is a morphism of Jordan triples and
$f(e_1,\ldots, e_r)$ is a Jordan frame.

\Proof. Let $x_0 := f(0)$. Then $\ell(x) := f(x) - x_0$ defines an
isometric linear
map $\ell \: \oline\Delta^r \to V$. Since $\ell$ is linear and
isometric, it maps the
open unit ball $\Delta^r$ in $\C^r$ into the open unit ball ${\cal
D}$ of $(V,|\cdot|)$,
so that it also maps $\oline\Delta^r$ isometrically into $\oline{\cal D}$.

Let $f_1, \ldots, f_r$ denote the images of the canonical basis in
$\C^r$ under $\ell$.
Then the coordinate projections
$$\chi_j \: L := \span \{f_1, \ldots, f_r\} = \im(\ell) \to \C, \quad
\sum_j \lambda_j f_j \mapsto \lambda_j $$
are linear maps with $\|\chi_j\| = 1$ because $\ell \: \C^r \to L$ is
an isometric
inclusion. Using the Hahn--Banach Theorem, we find extensions
$\chi_j \: V\to \C$ with the same norm. Then the map
$$ \chi := (\chi_1, \ldots, \chi_r) \: V \to \C^r $$
satisfies $\|\chi\| = 1$ and $\chi \circ \ell =\id$. It follows in particular
that $\chi({\cal D}) \subeq \Delta^r$.

Since $\chi$ maps $\oline{\cal D}$ into $\oline\Delta^r$, we have an
order-preserving map
$$ \chi^* \: {\cal F}(\oline\Delta^r) \to {\cal F}(\oline{\cal D}), \quad
F \mapsto \chi^{-1}(F) $$
and the corresponding map
$$ \ell^* \: {\cal F}(\oline{\cal D}) \to{\cal F}(\oline\Delta^r), \quad
F \mapsto \ell^{-1}(F) $$
satisfies
$$ \ell^* \circ \chi^* = (\chi \circ \ell)^* = \id. $$
We conclude that $\chi^*$ is an order preserving injection. This
entails in particular, that
for each strictly increasing chain
$$ F_0  \subset F_1 \subset F_2 \subset \ldots \subset F_r $$
of faces of $\oline\Delta^r$, the images under $\chi^*$ form a
strictly increasing chain of
faces of $\oline{\cal D}$. Since $r$ is the rank of ${\cal D}$, the
maximal chains in
${\cal F}(\oline{\cal D})$ are of length $r$, which implies that
$\chi^*$ preserves the rank of faces. Since the rank of an element $x
\in \oline{\cal D}$
coincides with the rank of the face it generates, we further see that for
$z \in \oline\Delta^r$ we have
$$ \rk \ell(z) = \rk \Face(\ell(z)) = \rk \ell^*(\Face(z)) =
\rk(\Face(z)) = \rk z. $$
Therefore $\ell$ preserves the rank.

Moreover, $\ell$ maps the Shilov boundary $\T^r$, consisting of the
elements of maximal rank,
into the Shilov boundary $S$ of ${\cal D}$. The relation
$$ f(\oline\Delta^r) = x_0 + \ell(\oline\Delta^r) \subeq \oline{\cal D} $$
implies
$$ - x_0 + \ell(\oline\Delta^r)  =  -(x_0 + \ell(\oline\Delta^r))
\subeq \oline{\cal D}, $$
so that for each $z \in \T^r$ we have
$$ \ell(z) = \half((\ell(z)+ x_0) + (\ell(z)-x_0)) \in S, $$
so that $S = \Ext(\oline{\cal D})$ implies $x_0= 0$, and hence $f =
\ell$ is linear.

For $i \in \{1,\ldots, r\}$ we consider the corresponding face
$$ F := \{ z \in \oline\Delta^r \: z_i = 1\} \in {\cal F}(\oline\Delta^r). $$
Then $F$ is the closure of an $(r-1)$-dimensional affine polydisc, and
$f\res_{F} \: F \to \oline{\cal D}$ is an affine isometry into a face
$F_c \in {\cal F}(\oline{\cal D})$, where $c$ is a primitive tripotent
(Theorem~I.5, Prop.~II.3). Applying the first part of the
proof with ${\cal D}$ replaced by $\algint(F')$ to the corresponding map
$$\oline\Delta^{r-1} \to F_c - c, \quad z \mapsto
f(z_1, \ldots, z_{i-1}, 1, z_i, \ldots, z_r) - c,$$
we see that this map is linear, hence maps $0$ to $0$, which leads to
$f(e_i) = c$. For $i \not=j$ the element
$e_i + e_j \in\oline\Delta^r$ is contained in the face generated by
$e_i$, which implies
that $f(e_i + e_j) = f(e_i) + f(e_j)$ is contained in the face generated by
$f(e_i)$. From Theorem~I.5 we now derive
$$ f(e_j) = f(e_i + e_j) - f(e_i) \in V_0(f(e_i)), $$
so that the primitive tripotents $f(e_i)$, $i =1,\ldots, r$, are
mutually orthogonal.
Hence the linear map $f \: \C^r \to V$ is a morphism of Lie triples systems.
\qed

\Corollary III.6. Suppose that ${\cal D}_1 \subeq V_1$ and
${\cal D}_2 \subeq V_2$ are circular bounded symmetric domains of the
same rank.
Then any affine isometric map $f \: V_1 \to V_2$ mapping
$\oline{\cal D}_1$ into $\oline{\cal D}_2$ is linear and rank-preserving.

\Proof. Let $r := \rk {\cal D}_1 = \rk {\cal D}_2$ and fix a polycylinder
${\cal D}_0 := \Delta^r \subeq {\cal D}_1$ defined by a Jordan frame
$(c_1, \ldots, c_r)$. For $V_0 := \span \{c_1,\ldots, c_r\}$ we then
obtain by restriction an
isometric map $f_0 \: V_0 \to V_2$ mapping $\oline{\cal D}_0 \to
\oline{\cal D}_2$.
In view of Theorem III.5, this map is linear, which implies
$f(0) = f_0(0)= 0$, and thus $f$ is linear.

Moreover, $f_0$ is rank-preserving by Theorem~III.5, which implies
that $f$ is also
rank-preserving.
\qed

\Corollary III.7. If $r = \rank {\cal D}$, then any
isometric linear embedding $f \: \Delta^r \into {\cal D}$ is
equivariant in the sense that
there exists a subgroup $G_1 \subeq \Aut({\cal D}_0)$ and a
surjective homomorphism
$G_1 \to \Aut(\Delta^r)_0 \cong \PSU_{1,1}(\C)^r$ such that
$f$ is equivariant with respect to the action of $G_1$ on $\Delta^r$
and ${\cal D}$.

\Proof. If $(e_1,\ldots, e_r)$ is the canonical basis in $\C^r$, then
$(c_1,\ldots, c_r) := (f(e_1), \ldots, f(e_r))$ is a Jordan frame,
so that
$$ \g_1 := \sum_{j=1}^r \g_{c_j} \subeq \g $$
is isomorphic to $\su_{1,1}(\C)^r \cong \sL_2(\R)^r$ (see (I.8)), the
Lie algebra
of the group $\Aut(\Delta^r)_0 \cong \PSU_{1,1}(\C)$. We may now put
$G_1 := \la \exp \g_1 \ra \subeq G$, and the assertion follows.
\qed

\sectionheadline{IV. The Maslov index}

To define the integers classifying the $G$-orbits in $S\times S\times
S$, we need in particular the
Maslov index, a certain $G$-invariant function $\iota \: S\times
S\times S \to \Z$. In this section
we explain how the Maslov index can be defined for bounded symmetric
domains of tube type
which are not necessarily irreducible, hence extending the definition
given in [C\O{}01], [C\O{}03], [Cl04b]. Using Theorem~III.1, we
further derive a
list of properties of the Maslov index and show that it can be
characterized in an axiomatic
fashion by these properties. Actually this was our original
motivation to prove Theorem~III.1.

Let us first consider the case of the unit disc $\Delta$. Then the group $G$ is
$\PSU_{1,1}(\C)$ acting by homographies on $\Delta$, and its Shilov
boundary is the
unit
   circle $\T$. The {\it Maslov index\/}
$$\iota =\iota_{\T} : \T\times \T\times
\T\longrightarrow \Bbb Z$$
is defined by
\smallskip
$\bullet \quad\iota (x,y,z)=0$ if two of the elements of the triplet coincide.
\smallskip
$\bullet \quad\iota(x,y,z) = \pm1$ if $(x,y,z)$ is conjugate under $G$ to $
(1,-1,\mp i).$

If $\Delta^r$ denotes the $r$-polydisc, then the identity component of
$\Aut(\Delta^r)$ is $G=\PSU_{1,1}(\C)^r$
and the Shilov boundary of $\Delta^r$ is $\T^r$. The Maslov index $\iota =
\iota_{\T^r} : \T^r\longrightarrow \Bbb R$ is defined by
$$ \iota((x_1,x_2,\dots,x_r), (y_1,y_2,\dots,y_r), (z_1,z_2,\dots,z_r)):=
\iota(x_1,y_1,z_1)+\iota(x_2,y_2,z_2)+\dots+\iota(x_r,y_r,z_r)\ .$$

Now consider an irreducible bounded symmetric domain $\cal D$ of tube type with
Shilov boundary~$S$. The Maslov index $\iota=\iota_S : S\times S\times
S\longrightarrow \Bbb Z$ is defined in [C\O{}01], [C\O{}03], [Cl04b]. As the
definition is involved, we won't repeat it here, but it has the following
property, which, in the light of
Theorem III.1 and because of the invariance of this index under $G$,
is characteristic: For any Jordan frame
$(c_1,c_2,\dots, c_r)$, let
$$T=\Big\{\sum_{j=1}^r t_j c_j \: |t_j| = 1, 1\leq j\leq r\Big\}$$
be the $r$-torus which is the Shilov boundary of the associated $r$-polydisc.
Then for any three points $x,y,z$ in $T$, one has
$$ \iota_S(x,y,z)= \iota_T(x,y,z). \leqno(4.1) $$

Last, we extend now the definition of the Maslov index to any bounded symmetric
domain $\cal D$ in the following way. Assume that $ {\cal D} = {\cal D}_1\times
{\cal D}_2\times \dots \times{\cal D}_m$ is the decomposition of $\cal D$ as a
product of irreducible domains. Then the identity component of the group of
biholomorphic automorphisms  of $\cal D$ is the product
$$ G= \Aut({\cal D}_1)_0\times \Aut({\cal D}_2)_0\times \dots \times
\Aut({\cal D}_m)_0,$$
and the Shilov boundary $S$ of $\cal D$ is the product
$S=S_1\times S_2\times \dots \times S_m$
of the corresponding Shilov boundaries. Then the Maslov index $\iota =
\iota_S$ is defined by
$$\iota(x,y,z) :=
\iota_{S_1}(x_1,y_1,z_1)+\iota_{S_2}(x_2,y_2,z_2)+\dots+\iota_{S_r}(x_l,y_l,z_l)\
.$$

\Theorem IV.1. The Maslov index has the following properties :
\litemindent=0.9cm
\litem{(M1)} It is invariant under the group $G$.
\litem{(M2)} It is an alternating function with respect to any
permutation of the
three arguments.
\litem{(M3)} It satisfies the cocycle property
$\iota(x,y,z) = \iota(x,y,w)- \iota(x,z,w) +\iota(y,z,w)$.
\litem{(M4)} It is additive in the sense that if ${\cal D} = {\cal
D}_1 \times {\cal
D}_2$, so that $S=S_1\times S_2$, then
$$ \iota_S(x,y,z) =
\iota_S((x_1,x_2),(y_1,y_2),(z_1,z_2))=\iota_{S_1}(x_1,y_1,z_1)+\iota_{S_2}(x_2,y_2,z_2)\
.$$
\litem{(M5)} If $\Phi: {\cal D}_1\longrightarrow {\cal D}_2$ is an
equivariant holomorphic
embedding of bounded symmetric domains of tube type of equal rank, then
$ \iota_{S_2}\circ\Phi = \iota_{S_1}$.
\litem{(M6)} It is normalized by $\iota_\T(1,-1,-i) = 1$ for the
Shilov boundary $\T$
of the unit disc $\Delta$.
\litemindent=0.7cm

\Proof. Properties (M1)-(M3) are known for irreducible domains
([C\O{}01], [Cl04]), and the extension
of these properties to products of irreducible domains is obvious.
Property (M4) obviously holds by the way we have defined the Maslov index.

For Property (M5), let $r$ be the common rank of the two domains.
We may assume that ${\cal D}_1$ and ${\cal D}_2$ are given in a
circular realization
as unit balls in spaces $V_1$, resp., $V_2$. Then
$\phi(0) \in {\cal D}_2$, and there is some $g_2 \in G_2 :=
\Aut({\cal D}_2)_0$ with
$g_2.\phi(0) = 0$. Then $\psi(z) := g_2.\phi(z)$ defines an
equivariant embedding
${\cal D}_1 \to {\cal D}_2$ which is linear because $\psi(0) = 0$.

Let $(x,y,z) \in S_1$ and pick $g_1 \in G_1 := \Aut({\cal D}_1)_0$ such that
$g_1.(x,y,z)$ is contained in the span of a Jordan frame $(c_1, \ldots, c_r)$
(Theorem~III.1), hence in the Shilov boundary $T_1$ of the corresponding
polydisc $\Delta^r$ in ${\cal D}_1$. From the equivariance of $\phi$
we derive the
existence of some $\tilde g_1 \in G_2$ with $\phi \circ g_1 = \tilde
g_1 \circ \phi$.
Then $\psi(\Delta^r)$ is a maximal polydisc in
${\cal D}_2$ with Shilov boundary $T_2 := \psi(T_1)$, so that (4.2)
implies that
$$ \eqalign{ \iota_{S_1}(x,y,z)
&= \iota_{S_1}(g_1.x,g_1.y,g_1.z)
= \iota_{T_1}(g_1.x,g_1.y,g_1.z) \cr
&= \iota_{T_2}(\psi(g_1.x),\psi(g_1.y),\psi(g_1.z))
= \iota_{S_2}(\psi(g_1.x),\psi(g_1.y),\psi(g_1.z)) \cr
&= \iota_{S_2}(g_2\phi(g_1.x),g_2\phi(g_1.y),g_2\phi(g_1.z))
= \iota_{S_2}(\phi(g_1.x),\phi(g_1.y),\phi(g_1.z)) \cr
&= \iota_{S_2}(\tilde g_1\phi(x),\tilde g_1\phi(y),\tilde g_1 \phi(z))
= \iota_{S_2}(\phi(x),\phi(y),\phi(z)). \cr} $$

Property (M6) is a consequence of the definition.
\qed

\Remark IV.2. Note that (M2) and (M3) mean that $\iota_S$ is a $\Z$-valued
Alexander--Spanier $2$-cocycle on $S$.
\qed

Before we turn to the general case in the following section,
we recall the classification of triples
in the circle, the Shilov boundary of the unit disc:

\Example IV.3. We consider the case $\Delta := \{ z \in \C \: |z| < 1\}$.
Then $G = \PSU_{1,1}(\C)$ acts by
$$ \Big[\pmatrix{a & b \cr c & d \cr}\Big].z = (az + b)(cz+d)^{-1}. $$
The Shilov boundary is $S = \T= \{ z \in \C \: |z| = 1\}$.
Identifying $S$ with the projective line $\P_1(\R)$ and $G$ with $\PSL_2(\R)$,
we immediately see that there are exactly two $G$-orbits in $S \times
S$, represented by
$$ (1,1) \quad \hbox{ and } \quad (1,-1), $$
i.e., the diagonal in $S \times S$ and the set $(S \times S)_\top$ of
transversal pairs.
Since the action of $G$ on $S$ preserves the orientation of a triple, it
follows that we have $6$ orbits in $S\times S\times S$, represented by
$$ (1,1,1), \quad
(1,1,-1), \quad (1,-1,1), \quad (1,-1,-1), \quad (1,-1,-i)
\quad \hbox{ and } \quad (1,-1,i).
\qeddis

\Remark IV.4. As a function assigning to any triple in the Shilov
boundary of any
bounded symmetric domain ${\cal D}$ an integer, the Maslov index is uniquely
determined by the properties (M1), (M2) and (M4)-(M6).

In view of Example~IV.3, the Maslov index for ${\cal D} = \Delta$
is uniquely determined by (M1), (M2) and (M6). By (M4) it is also determined
for polydiscs.

If ${\cal D}$ is any bounded symmetric domain of rank $r$ and $(s_1,
s_2, s_3) \in S\times S\times S$,
then Theorem~III.1 implies that it can be conjugate by some $g \in G$
to a triple in
the Shilov boundary $T \cong \T^r$ of a maximal polydisc, so that
Corollary~III.7,
(M1) and (M5) lead to
$$ \iota_S(s_1, s_2, s_3) = \iota_S(g.s_1, g.s_2, g.s_3) =
\iota_T(g.s_1, g.s_2, g.s_3). $$
We conclude that $\iota_S$ is determined uniquely by (M1), (M2),
together with (M3)-(M6).
\qed

\subheadline{A classical case: the Lagrangian manifold}

Let $E$ be a real vector space of dimension $2r$ and $\omega$ be a
symplectic form
on $E$. The symplectic group $\Sp(E,\omega)$
is the group of linear automorphisms which preserve
$\omega$. A {\it Lagrangian\/} is a maximal totally isotropic subspace
of $E$, hence of dimension $r$.
The set $\Lambda_r$ of all Lagrangians is a compact submanifold of the
Grassmannian
$\Gr_r(E)$ of $r$-dimensional subspaces of $E$.
Then the group $G:=\PSp(E,\omega) := \Sp(E,\omega)/\{\pm \1\}$ acts
transitively and effectively
on $\Lambda_r$. Choosing a symplectic basis in
$E$, we may identify $E$ with $\R^r\times \R^r$, the symplectic form being
the standard one, namely
$$ \omega((\xi,\eta),(\xi',\eta'))=\xi^\top\eta'-\eta^\top\xi'. \leqno(4.2) $$

Let us consider the complex vector space $V=\Sym_r(\C)$ of complex $r\times r$
symmetric matrices, and let $\cal D$ be the unit ball with respect to
the operator norm.
The space $V$ is an involutive unital Jordan algebra with real form
$\Sym_r(\R)$, involution $z^* = \oline z$ and Jordan product
$x * y := \half(xy + yx)$. The spectral norm
on $V$ coincides with the operator norm, and the unit ball is then a
bounded symmetric
domain. To make connection with symplectic geometry, observe that the
graph of a
symmetric matrix is a complex isotropic subspace in $\C^r\times \C^r$ for the
symplectic structure (4.2). Let moreover $h$ be the
hermitian form on
$\C^r\times \C^r$ given by
$$ h((\xi,\eta),(\xi',\eta'))=\xi^\top\overline
{\xi'}-\eta^\top\overline {\eta'}
= (\xi')^* \xi - (\eta')^* \eta.$$
The hermitian form $h$ has signature $(r,r)$. Now to any $x\in V$,
associate its graph
$$ \ell_x=\{(\xi,x.\xi)\: \xi\in \C^r\}.$$
The condition that $x$ is in the unit ball is equivalent
to the fact that $ \1 - xx^*$ is positive definite, which in turn
implies that the
restriction of
$h$ to $\ell_x$ is positive definite. Conversely, any (complex)
Lagrangian in $\C^r\times \C^r$
on which the restriction of $h$ is positive definite is the graph of
some complex symmetric matrix in the unit ball. The Shilov boundary
of $\cal D$ is the
manifold of unitary symmetric matrices, and the corresponding graphs
are the (complex)
Lagrangians on which the restriction of the form $h$ is identically
$0$. Let $C$ be the
map from $\R^r\times \R^r$ to $\C^r\times \C^r$ given by
$$ C(\xi,\eta) = \Big({\xi+i\eta \over \sqrt 2},
{\xi-i\eta\over\sqrt2 }\Big)\ .$$
Then an elementary computation shows that the complexification of the image
under $C$ of a (real)
Lagrangian is a (complex) Lagrangian on which the restriction of $h$
is identically $0$,
and vice versa. This gives a one-to-one correspondence between
$\Lambda_r$ and $S$.
Moreover the natural action of $G$ on $\Lambda_r$ is transferred to
an action on $S$
and realizes an isomorphism of the real symplectic group and the
group $\Sp_{2r}(\C)\cap
\UU_{r,r}(\C)$, which generalizes the isomorphism of $\SL_2(\R)$ and
$\SU_{1,1}(\C)$.

The matrices $E_{11}, \ldots, E_{rr}$ form a
Jordan frame in $\Sym_r(\C)$. The corresponding $r$-torus is
$$ T := \Bigg\{\pmatrix
{e^{i\theta_1}&0&\dots&0\cr0&e^{i\theta_2}&\dots&0\cr\vdots&\vdots&\ddots&\vdots\cr0&0&\dots&e^{i\theta_r}}\: 
\theta_j\in \R, 1\leq j\leq r\Bigg\}\ .$$
The graph of an element of $T$ is the $r$-space generated by
$$ (e_1,e^{i\theta_1} e_1),  (e_2,e^{i\theta_2} e_2), \dots,
(e_r,e^{i\theta_r} e_r),$$
or equivalently by
$$ (e^{-i{\theta_1\over 2}}e_1, e^{i{\theta_1\over 2}}e_1),
(e^{-i{\theta_2\over 2}}e_2,
e^{i{\theta_2\over 2}}e_2),\dots,(e^{-i{\theta_r\over 2}}e_r,
e^{i{\theta_r\over 2}}e_r).$$
Observe that $ (e^{-i{\theta_j\over 2}}e_j, e^{i{\theta_j\over 2}}e_j) = C(\cos
{\theta_j\over 2} e_j, \sin{\theta_j\over 2}e_j)$ to get that the corresponding
Lagrangian  $\ell(\theta_1,\theta_2,\dots, \theta_r)$ in $\Lambda_r$ is
   generated by
$$ \Big(\cos {\theta_1\over 2}e_1, -\sin {\theta_1\over 2}e_1\Big),
\Big(\cos {\theta_2\over 2}e_2,-\sin {\theta_2\over 2}e_2\Big),
\ldots,\Big(\cos {\theta_r\over 2}e_r, -\sin {\theta_r\over 2}e_r\Big)\ .$$
In this case, one can then reformulate Theorem III.1  as follows.

\Theorem IV.5. Let $\ell_1,\ell_2,\ell_3$ be three arbitrary Lagrangians in a
symplectic vector space $E$ of dimension $2r$. Then there exists a
symplectic basis
$e_1,e_2,\dots, e_r,f_1,f_2,\dots,f_r$ such that each of the three
Lagrangians is
generated by
$$\cos \theta_1 e_1+\sin\theta_1 f_1,\cos \theta_2 e_2+\sin\theta_2
f_2,\dots,\cos
\theta_r e_r+\sin\theta_r f_r  $$
for appropriate choices of the $(\theta_j)_{1\leq j \leq r}$.
\qed

The classification result (Theorem~V.4 below) for the case
$S = \Lambda_r$ can also be found in [KS90, p.492].

\sectionheadline{V. The classification of triples}

In this section we complete the classification of $G$-orbits in the
set $S\times S\times S$ of triples
in $S$ by first assigning to each triples an increasing $5$-tuple of integers
$N = (n_1, n_2, n_3, n_4, n_5)\in \{0,\ldots, r\}^5$
depending only on its orbit. Then we exhibit for each such $5$-tuple
a standard triple
with this invariant, and finally we show that two different standard
triples belong to
different orbits.

\Definition V.1. To
any triple $(x_1,x_2,x_3)$ in $S\times S\times S$, we may
associate five integers:
\litem{(1)} the ranks of the three faces (cf.\ Remark~II.4):
$$n_{12}=\rank \Face(x_1,x_2),\quad
n_{2,3} =\rank \Face(x_2,x_3),\quad n_{3,1}= \rank
\Face(x_3,x_1)$$
\litem{(2)} the rank of the face generated by the
triple
$$n_{1,2,3} =\rank \Face(x_1,x_2,x_3)$$
\litem{(3)} the Maslov index $\iota(x_1,x_2,x_3)$.
\medskip
Clearly the action of $G$ preserves these integers.
\qed

When $x_1,x_2,x_3$ are contained in the boundary of a polydisc
(cf.\ Section III), then these integral invariants are easy to
compute (cf.\ Example~II.6).

\Lemma V.2. Let $e=\sum_{j+1}^r c_j$ be a Peirce decomposition of the
unit, and, for $\kappa = 1,2,3,$ let
$$ x_\kappa = \sum_{j=1}^r \xi_j^{(\kappa)} c_j, \quad
\hbox{ where } \quad
\vert \xi_j^{(\kappa)}\vert = 1 \quad \hbox{for all } \quad j \in
\{1, \ldots, r\}. $$
Then
$$ n_{\kappa,\kappa'} = |\{j \: \xi_j^{(\kappa)}= \xi_j^{(\kappa')}\}|, \quad
n_{1,2,3} = | \{ j \: \xi_j^{(1)} = \xi_j^{(2)}=\xi_j^{(3)}\}|, $$
and
$$ \iota(x_1,x_2,x_3) =\sum_{j=1}^r \iota
(\xi_j^{(1)},\, \xi_j^{(2)},\, \xi_j^{(3)}).
\qeddis

\Definition V.3. We now describe the {\it standard triples\/}
associated to a (fixed)
Jordan frame $(c_1, \ldots, c_r)$. Let
$N=(n_1,n_2,n_3,n_4,n_5)$ be a
$5$-tuple of integers such that
$$0\leq n_1\leq n_2\leq n_3\leq n_4\leq n_5\leq r \ .$$

Then the {\it standard triple of type\/} $N$ is the triple
$(x_1^N,x_2^N,
x_3^N)
$  defined by
$$ x^{N}_1 = e_r = c_1 + \ldots + c_r,\qquad
x^{N}_2= c_1+c_2+\dots +c_{n_2}-c_{n_2+1}-\dots  -c_r,$$
$$ x^N_3= c_1+\dots +c_{n_1}-c_{n_1+1}-\dots - c_{n_3}+c_{n_3+1}+\dots
+c_{n_4} -i c_{n_4+1}-\dots -ic_{n_5}+ic_{n_5+1}+\dots +ic_r\ .$$

For this triple, one has
$$ n_{1,2,3}= n_1, \quad
n_{1,2}= n_2,\quad n_{1,3}= n_1+n_4-n_3,\quad n_{2,3}= n_1+ n_3-n_2,$$
and
$$ \iota(x_1^N,x_2^N, x_3^N)=n_5-n_4-(r-n_5) = 2n_5 - n_4 - r.
\qeddis

\Theorem V.4. If ${\cal D}$ is an irreducible bounded symmetric
domain of tube type, then
any triple in $S$ is conjugate to one and only one of the standard triples.

\Proof. For the standard triples we have
$$ n_1 = n_{1,2,3}, \quad n_2 = n_{1,2}, \quad
n_3 = n_{2,3} + n_2 - n_1 = n_{2,3} + n_{1,2} - n_{1,2,3}, \leqno(5.1) $$
$$ n_4 = n_{1,3} + n_3 - n_1 = n_{1,3} + n_{2,3} + n_{1,2} - 2
n_{1,2,3}, \leqno(5.2) $$
and
$$ n_5
= \half(\iota(x_1^N,x_2^N, x_3^N) + n_4 + r)
= \half(\iota(x_1^N,x_2^N, x_3^N) + r + n_{1,3} + n_{2,3} + n_{1,2}
- 2 n_{1,2,3}). \leqno(5.3) $$
Since the numbers $n_{1,2,3}$, $n_{1,2}$, $n_{2,3}$, $n_{3,1}$ and
the Maslov index
are $G$-invariant, it follows that for different values of $N$,
the corresponding standard triples are not conjugate under $G$.

To show, conversely, that each triple $(e,f,h) \in S\times S\times S$
is conjugate to a standard triple, we first use
Theorem~III.1 to see that we may w.l.o.g.\ assume that
$(e,f,h)$ is contained in the torus
$$ T := \Big\{ \sum_{j=1}^r \lambda_j c_j \: (\forall j)\ |\lambda_j|
=1\Big\} $$
defined by the Jordan frame $(c_1, \ldots, c_r)$.
It is the Shilov boundary
of the polydisc
$$ \Delta^r := \Big\{ \sum_{j=1}^r \lambda_j c_j \: (\forall j)\
|\lambda_j| < 1\Big\}. $$

We write
$$ e = \sum_{j=1}^r \xi_j^e c_j, \quad
f = \sum_{j=1}^r \xi_j^f c_j \quad \hbox{ and } \quad
h = \sum_{j=1}^r \xi_j^h c_j. $$
  From I.8 it follows that every element of $\Aut(\Delta^r)_0 \cong
\PSU_{1,1}(\C)^r$
is the restriction of an element of $\Aut({\cal D})_0$, because
$$ \g_{c_1} + \ldots + \g_{c_r} \cong \su_{1,1}(\C)^r = \L(\Aut(\Delta^r)) $$
is a subalgebra of $\g = \L(G)$. We may therefore assume that
$\xi_j^e = 1$ for each $j$. Let
$$n_2 := |\{ j \: \xi^e_j = \xi^f_j\}|
=  |\{ j \: \xi^f_j = 1\}|. $$
Since each permutation of the set $\{c_1, \ldots, c_r\}$ is induced
by an element
of $K$, which acts transitively on the set of Jordan frames, we may
w.l.o.g.\ assume that
$$ f = c_1 + c_2 + \ldots + c_{n_2} - c_{n_2+1} - \ldots - c_r $$
because the $\Aut(\Delta)_0$-orbits in $\T \times \T$ are represented by
$(1,1)$ and $(1,-1)$ (Example~IV.3).

Let $n_1 := |\{ j \: \xi^e_j = \xi^f_j = \xi^h_j\}|$ and write
$$n_4 := |
\{ j \: \xi^e_j = \xi^f_j\ \hbox{or}\ \xi^e_j = \xi^h_j \ \hbox{or}\
\xi^f_j = \xi^h_j\}| $$
for the number of components in which at least
two elements of $\{e,f,h\}$ have the same entries. Then
$h$ has precisely $n_1$ entries $1$ among the first $n_2$,
and we may w.l.o.g.\ assume that they arise in position $j = 1,\ldots, n_1$.
We may likewise assume that the components of $e,f$ and $h$ are
mutually different for
$j > n_4$. Then the entries of $h$ in positions
$n_1 + 1, \ldots, n_2$ can be moved by elements of the group
$\Aut(\Delta)_0^{n_2 - n_1}$ acting on these components to~$-1$.
For $j \in \{n_2 + 1,\ldots, n_4\}$ the $j$-th component of $h$
equals either $1$ or $-1$. Moving
the $1$-entries with some element of $K_e$ permuting $\{c_1,\ldots, c_r\}$
to the rightmost positions, we get entries $-1$ for
$j = n_1 + 1,\ldots, n_3$ for some $n_3$ satisfying $n_2 \leq n_3 \leq n_4$.
For $j > n_4$ we then have $\Im \xi^h_j \not=0$, and after permuting the Jordan
frame, we may assume that for some $n_5 \geq n_4$ we have
$\Im \xi^h_j < 0$ for $j = n_4 + 1,\ldots, n_5$ and
$\Im \xi^h_j > 0$ for $j > n_5$. We finally use elements of $\Aut(\Delta)_0$
fixing $1$ and $-1$ to move each entry with negative imaginary part to $-i$ and
the others to $i$ (cf.~Example~IV.3).
This proves that each triple is conjugate to a standard triple.
\qed

\Remark V.5. In Theorem~V.4, we have classified the $G$-orbits in the
space of triples
in $S$ by the set of all $5$-tuples
$N = (n_1, n_2, n_3, n_4, n_5) \in \{0,\ldots, r\}$
satisfying the monotonicity condition
$$ n_1 \leq n_2 \leq n_3 \leq n_4 \leq n_5. $$
The description the standard triples shows that each such tuples arises via
(5.1)-(5.3).
We claim that for the $5$-tuple
$$ (r_0, r_1, r_2, r_3, d) := \big(n_{1,2,3}, n_{1,2}, n_{2,3}, n_{3,1},
\iota(x_1^N,x_2^N,x_3^N)\big) $$
of integers we then have
\litem{(P1)} $0 \leq r_0 \leq r_1, r_2, r_3 \leq r$.
\litem{(P2)} $r_1 + r_2 + r_3 \leq r + 2 r_0$.
\litem{(P3)} $|d| \leq r + 2 r_0 - (r_1 + r_2 + r_3)$.
\litem{(P4)} $d \equiv r + r_1 + r_2 + r_3 \mod 2$.

In fact, (P1) is clear,
$$ r_1 + r_2 + r_3 = n_4 + 2 r_0 \leq r + 2 r_0, $$
$$ |d| = |n_5 - n_4 -(r - n_5)| \leq n_5 - n_4 + r - n_5 = r - n_4 =
r + 2 r_0 - r_1 - r_2 - r_3, $$
and
$$ d= n_5 - n_4 - (r - n_5) \equiv n_4 + r \equiv r + r_1 + r_2 +
r_3\ \mod 2. $$

Suppose, conversely, that $(r_0, r_1, r_2, r_3, d) \in \Z^5$
satisfies (P1)-(P4).
We then define
$$ n_1 := r_0,\quad n_2 := r_1,\quad n_3 := r_2 + r_1 - r_0, \quad
n_4 := r_3 +  r_2 + r_1 - 2 r_0 $$
and
$$ n_5 = \half(d+ r_3 +  r_2 + r_1 + r) -r_0. $$
Then
(P4) implies $n_5 \in \Z$.
  From (P1/2) we immediately get
$0 \leq n_1  \leq n_2 \leq n_3 \leq n_4 \leq r.$
Further (P3) leads to $|d| \leq r - n_4$, and
$n_4 \leq n_5$ follows from
$$ 2 n_5 = d + r_3 + r_2 + r_1 + r - 2 r_0 = d + r + n_4 \geq r + n_4
- (r - n_4) = 2 n_4. $$
This is turn implies
$n_5 = \half(r + d + n_4) \leq r.$
\qed

The conditions (P1)-(P4) are well known conditions describing the
classification of
triples of Lagrangian subspace of symplectic vector spaces ([KS90]).

\sectionheadline{VI. Classification of orbits in $S\times S$}

In this section we describe how the classification of $G$-orbits in
$S \times S$
can be derived from the Bruhat decomposition of $G$, resp., the
description of the
orbits of the maximal parabolic subgroup $G_e$ in $G$ with $G/G_e \cong S$.

Throughout this section we assume $\cal D$ to be irreducible.
Let $(c_1,c_2,\dots, c_r)$ be a Jordan frame and put
$$ \varepsilon_k=c_1+c_2\dots+c_k-c_{k+1}-\dots-c_r \quad
\hbox{ for } \quad k = 0,\ldots, r. $$
Moreover let $e= c_1+\dots +c_r=\varepsilon_r$, and observe that $
\varepsilon_0 = -e$. The vector space
$$ \goth a = \bigoplus_{j=1}^r \Bbb Rc_j$$ is a maximal flat in $V$
in the sense of
Loos ([Lo77]) and can be thought of as a  Cartan subspace in
the tangent space of $\cal D$ at the origin. The corresponding vector fields
form a Cartan subspace of $\goth p$. Denoting by
$\gamma_j$ the
$j$-th coordinate in
$\goth a$ with respect to the basis $(c_1,c_2,\dots, c_r)$, it is known that
the (restricted) roots of $ (\goth g,\goth a)$ are $ \pm \gamma_j\pm
\gamma_k , \pm
2\gamma_j, 1\leq j\neq k\leq r$ and, in addition,
$\pm\gamma_j, 1\leq j\leq r$ in the non tube type case. We
choose as positive Weyl chamber in $\goth a$ the one defined by the
inequalities
$$ \gamma_1\geq\gamma_2\geq \dots \geq \gamma_r\geq 0,$$
so that the corresponding simple roots are
$$\gamma_1-\gamma_2,\ \gamma_2-\gamma_3,\dots, \gamma_{r-1}-\gamma_r,\
   \gamma_r\ .$$
The Weyl group $W$ is isomorphic to the semi-direct product $\goth S_r\ltimes
\Bbb Z_2^r$, where $\goth S_r$ acts by permutation of the coordinates
$\gamma_j$, and the $j$-th factor $\Bbb Z_2$ acts by changing the sign of the
$j$-th coordinate.

The stabilizer $G_e$
of the point $e\in S$ is known to be a maximal parabolic subgroup
(cf.~Sect.~I).
It is the
standard parabolic subgroup associated to the subset
$$\Theta =\{\gamma_1-\gamma_2,\
\gamma_2-\gamma_3,\dots, \gamma_{r-1}-\gamma_r\}$$
of the set of simple roots. The subgroup $W^{\Theta}$of
$W$ generated by the reflections associated to the roots in $\Theta$
is just
$\goth S_r$, and  double cosets in $W^\Theta \backslash W/W^\Theta$
correspond to orbits of $\goth S_r$ in $\Z_2^r$, which
are
characterized by their number of sign changes. In particular, this
shows that the
elements $ \varepsilon_j, 0\leq j\leq r$, form a set of representatives of the
$W^\Theta$-orbits in $W.e$.

\Theorem VI.1. There are $r+1$ orbits of $G$ in $S\times S$. A set of
representatives of these orbits is given by the pairs
$ (e,\varepsilon_k), 0\leq k\leq r$.

\Proof. As $G$ acts transitively on $S$, any orbit of $G$ in $S\times
S$ meets the subset $\{e\}\times S$. So the statement amounts to show that a
$G_e$-orbit in $S$ contains  $\varepsilon_k$ for some $k,0\leq k\leq
r$. By Bruhat's
theory, the orbits of the parabolic subgroup $G_e$ of $G$ are in
one-to-one correspondence
with the $W^\Theta$-double cosets in $W$. In view of the preceding
discussion, this
shows the result.
\qed

\Remark VI.2. The open orbit in $S$ under the $G_e$-action (the big
Bruhat's cell)
corresponds to the point $-e$ and is nothing but the set of all
points in $S$ transversal to $e$.
\qed

\Definition VI.3. For $(x,y)\in S\times S$ we define their {\it
transversality index\/}
$\mu(x,y)$ to be the unique number $k \in \{0,\ldots, r\}$ such that
$(x,y)$ belongs to the
$G$ orbit of
$(e,\varepsilon_k)$. Clearly, the transversality index is invariant by the
action of $G$, and two pairs are conjugate if and only if they have the same
transversality index. Moreover, a pair $(x,y)$ is transversal if and
only if its
transversality index is $0$.
\qed

\Theorem VI.4.  A pair $(x,y)\in S\times S$ has transversality index $k$ if
and only if the face $F(x,y)$ generated by $x$ and $y$ has rank $k$.

\Proof. For $0\leq k\leq r$ let $e_k=c_1+c_2+\dots + c_k$. Then the
face generated by $e$ and $\varepsilon_k$ is
$$ \Face(e,\varepsilon_k) = (e_k+V_0(e_k)) \cap \overline {\cal D}, $$
which has rank $k$. As any pair in $S\times S$ is conjugate to one of the pairs
$(e,\varepsilon_k)$, the theorem follows immediately.
\qed

\sectionheadline{ Appendix: Bounded symmetric domains and tube type domains}

In this appendix, we briefly review the relation between {\it bounded symmetric
domains} and {\it positive hermitian Jordan triple systems} on one
hand, and the
relation between {\it bounded symmetric domains of tube type} and
{\it euclidean
Jordan algebras} on the other hand. Main references are [Lo77]
for (hermitian) Jordan triples and [FK94] for (euclidean) Jordan algebras.
\msk

A {\it hermitian Jordan
triple system\/}
$V$ is a finite dimensional complex vector space, together with a map
$\{\cdot,\cdot,\cdot\} : V\times V\times V\longrightarrow V $, such that
$\{x,y,z\} $ is
complex linear in $x$ and $z$, conjugate-linear in $y$, and such that
$$ \{x,y,z\} = \{z,y,x\} \leqno (JT1)$$
$$ \{a,b,\{x,y,z\}\} = \{\{a,b,x\},y,z\}-\{x,\{b,a,y\},z\}+
\{x,y,\{a,b,z\}\} \leqno (JT2)$$
for all $a,b,x,y,z\in V$.

\medskip
For $x,y\in V$ denote by $x\square y$ the linear endomorphism of $V$
defined by
$$ (x\square y)\,z = \{x,y,z\}$$
and by $Q(x)$ the conjugate linear endomorphism of $V$ defined by
$Q(x)z= \{x,z,x\}.$
Define the {\it trace form\/} $B$ on $V$ by
$B(x,y) = \tr (x\square y).$
The Jordan triple system $V$ is said to be {\it non degenerate\/} if, as a
sesquilinear form, $B$ is non degenerate. If this is the
case, then $B$ is hermitian (i.e. $ B(x,y) = \overline {B (y,x)}$ for all
$x,y\in V$). If moreover $B$ is {\it positive definite}, then
$V$ is said to be a {\it positive hermitian Jordan triple system}.
\medskip
Let $V$ be a positive hermitian Jordan triple system. An element
$c\in
V$ is said to be a {\it tripotent\/} if $\{c,c,c\}=c$.
For a tripotent $e \in V$ let
$V_j := V_j(e)$ denote the $j$-eigenspace of the operator $2 e \square e$.
Then we obtain the corresponding {\it Peirce decomposition of $V$}:
$$ V  = V_0 \oplus V_1 \oplus V_2 $$
([Lo77, Th.~3.13]).
\medskip
There is a (partial) order relation on tripotents. For  two
tripotents $c,d \in V$, we define $c \prec d$
if there exists a tripotent $c'$, such that
\litem{(i)} $c\square c'=0$ (orthogonality of $c$ and $c'$)
\litem{(ii)} $d=c+c'$.

A non zero tripotent is said to be {\it primitive\/} if it is minimal
among non zero tripotents for this order. Any tripotent
$c$ can be written as a sum of parwise orthogonal primitive tripotents,
say
$c= c_1+c_2+\dots +c_k$. The number $k$ of primitive tripotents in such
a decomposition of $c$ depends only on $c$ and is called the {\it rank of
$c$}.
\medskip
A {\it Jordan frame} of $V$ is a maximal family $(c_1,c_2,\dots,
c_r)$ of orthogonal
primitive tripotents. All Jordan frames have the same number of
elements called the
{\it rank of $V$}. For any Jordan frame $(c_1,c_2,\dots, c_r)$, the sum $
e=\sum_{j=1}^r c_j$ is a maximal tripotent of $V$, and all maximal
tripotents are
obtained this way.
\medskip
One of the main results in the theory of positive hermitian Jordan
triple system is
the {\it spectral theorem}.

\Proposition A.1. For any $x\in V$, there exists a Jordan frame
$(c_1,c_2,\dots, c_r)$ and positive real numbers $\lambda_j, 1\leq j
\leq r$, such
that
$$ x=\sum_{j=1}^r \lambda_jc_j\ .\leqno (A.1)$$
The $\lambda_j$ are unique up to a permutation.
\rm

The identity $(A.1)$ is called a {\it spectral decomposition\/} of $x$. The
$\lambda_j$ are called the {\it eigenvalues\/} of
$x$. The largest eigenvalue is the {\it spectral norm\/} of $x$,
denoted by $|x|$. As notation suggests, the map $x\mapsto |x|$ is
a complex Banach norm on $V$.

\Theorem A.2. The unit ball of $(V, |\cdot|)$ is a
bounded symmmetric domain. Conversely, any bounded symmetric domain is
holomorphically equivalent to such a unit ball.

\rm There is a subclass of symmetric bounded domains, the domains of 
{\it tube type}.
They are associated to a subclass of positive hermitian Jordan triple systems,
obtained by complexification from {\it euclidean Jordan algebras}.
\medskip
A {\it euclidean Jordan algebra} $E$  is a real finite dimensional
euclidean vector space $E$ with an inner product $\la \cdot, \cdot
\ra$, a bilinear
map $E\times E\longrightarrow E$ and an element $e\in E$ such that
$$ xy = yx, \qquad  ex=x, \quad x^2(xy) = x(x^2y) \quad \hbox{ and } \quad
\la xy,z\ra = \la y,xz\ra$$
for all $x,y,z\in E$.
Let $V=E^\Bbb C$ be the complexification of $E$, and extend the 
Jordan product from $E$
in a $\Bbb C$-bilinear way to $V$. Denote by $z\mapsto
\overline z
$ the conjugation of $V$ with respect to $E$.  For $x,y,z\in V$, let
$$ \{x,y,z\} := (x\overline y)z+x(\overline y z)-\overline y (xz)\ .\leqno
(A.3)$$ This endows $V$ with a structure of positive hermitian Jordan triple
system. The element $e$ is a tripotent of $V$. It satisfies $e\square e
= \id_E$, so that $V_0(e)=\{0\}$ (hence $e$ is a maximal tripotent),
$V_1(e)=\{0\}$ and $ V=V_2(e)$.
\msk
Among positive hermitian Jordan triple sytems, those coming from
euclidean Jordan
algebras are characterized by this last property. Let
$V$ be a positive hermitian Jordan triple system, and let $e$ be a
maximal tripotent. By
maximality of $e$,
$V_0(e)=\{0\}$. Assume further that $V_1(e)=\{0\}$, so that $
V=V_2(e)$. Now $Q(e)$
is a conjugate linear involution of
$V$. Its fixed points set $E=\{x\in V : Q(x) = x\}$ is a real vector space. For
$x,y\in E$, define
$$ xy = \{x,e,y\}\ .\leqno(A.4)$$
With the product defined by $(A.4)$ and the inner product induced by
$B$, $E$ is
then a euclidean Jordan algebra,
$V$ is the complexification of $E$ and the Jordan triple product on $V$ can be
recovered by formula (A.3) from the Jordan algebra product on $E$.
\medskip
An element $c\in E$ is called an idempotent if $c^2=c$. A Jordan frame in
$E$ is a maximal
set of orthogonal minimal idempotents. The number of elements in a
Jordan frame is
equal to $r$, the rank of the Jordan algebra $E$, and if
$(c_1,c_2,\dots, c_r)$ is a
Jordan frame, then $e=c_1+c_2+\dots +c_r$. A tripotent $c$ for the
associated triple
Jordan system structure on $V$ is of the form $ c=\sum_{j=1}^r
\lambda_j c_j$, for a
certain Jordan frame $(c_1,c_2,\dots, c_r)$ of $E$ and for each
$j,1\leq j\leq r$,
$|\lambda_j|=1$ or
$\lambda_j=0$. A maximal tripotent $z$ of $V$ is of the form with $
x=\sum_{j=1}^r
\lambda_j c_j$, with $|\lambda_j|=1$ for all
$j, 1\leq j\leq r$, so that $z$ as an element of the complex Jordan
algebra $V$ is
invertible, and satisfies $\overline z=z^{-1}$.

The
corresponding bounded symmetric domain is described as before by
$$ {\cal D} = \{z\in V : |z|<1\}\ .$$
The domain $\cal D$ can be shown to be holomorphically equivalent to
a tube domain.
If $E^+$ is the interior of the {\it cone of squares} of $E$, then by
the Cayley
transform $C_e$, the domain $\cal D$ is mapped to
$$ {\cal D}^C = C_e({\cal D}) = \{v\in V,\ \Re (v)>0\}= E^+\oplus iE$$
The domain $ {\cal D}^C$ is a tube domain in $V$, which is the
justification for
calling
$\cal D$ a bounded symmetric domain of tube type.
\msk
\nin The description of maximal tripotents of $V$ we gave supra shows
that the Shilov
boundary can be described as
$$ S= \{ z\in V : \overline z= z^{-1}\}\ .\leqno(A.5) $$
Hence the Shilov boundary $S$ is a totally real submanifold of $V$
with $\dim_\Bbb R
S =
\dim_\Bbb C V$.
This last condition is another characterization of bounded symmetric
domains of tube
type inside the family of bounded symmetric domains. In fact if $\cal D$ is a
bounded symmetric domain, then its Shilov boundary $S$ is a real
submanifold of $V$,
and its dimension satisfies $$\dim_\Bbb R S \geq {1\over 2} \dim_\Bbb
R V=\dim_\Bbb C V\ .$$
Equality is obtained if and only if $\cal D$ is of tube type.

\def\entries{

\[Be00 Bertram, W., ``The Geometry of Jordan and Lie Structures,''
Lecture Notes in Math. {\bf 1754}, Springer-Verlag, 2000

\[BN05 Bertram, W., and K.-H. Neeb, {\it Projective completions of
Jordan pairs,
Part II: Manifold structures and symmetric spaces}, Geom. Dedicata
{\bf 112:1} (2005), 75--115

\[BS64 Borel, A., and J.-P. Serre, {\it Th\'eor\`emes de finitude en 
cohomologie galoisienne}, Comment. Math. Helv. {\bf 39} (1964), 
111--164

\[CLM94 Cappel, S.E., R. Lee, and E.Y. Miller, {\it On the Maslov Index},
Comm. Pure Appl. Math. {\bf 47} (1994), 121--186

\[Cl04a Clerc, J.-L., {\it The Maslov triple index on the Shilov
boundary of a classical
domain}, J. Geom. Physics {\bf 49:1} (2004), 21--51

\[Cl04b ---, {\it L'indice the Maslov g\'en\'eralis\'e},
Journal de Math. Pure et Appl. {\bf 83} (2004), 99--114

\[Cl05 ---, {\it An invariant for triples in the Shilov
boundary of a bounded symmetric domain}, submitted

\[C\O{}01 Clerc, J-L., and B. \O{}rsted, {\it The Maslov index
revisited}, Transformation
Groups {\bf 6} (2001), 303--320

\[C\O{}03 ---, {\it The Gromov norm of the Kaehler class and the
Maslov index}, Asian J. Math. {\bf 7} (2003), 269--296

\[FMS04 Falbel, E., Marco, J.-P., and F. Schaffhauser, {\it
Classifying triples of
Lagrangians in a hermitian vector space}, Topology Appl. {\bf 144}
(2004), 1--27

\[FK94 Faraut, J., and A. Koranyi, ``Analysis on Symmetric Cones",
Oxford Mathematical Monographs, Oxford University Press, 1994

\[HR99 Hille, L., and G.~R\"ohrle, {\it A classification of parabolic
subgroups of classical
groups with a finite number of orbits on the unipotent radical},
Transf. Groups {\bf 4:1} (1999), 35--52


\[KS90 Kashiwara, M., and P. Schapira, ``Sheaves on Manifolds,''
Grundlehren der math. Wiss.
{\bf 292}, Springer-Verlag, 1990

\[Li94 Littelmann, P., {\it On spherical double cones}, J. Algebra
{\bf 166} (1994), 142--157

\[Lo77 Loos, O., ``Bounded Symmetric Domains and Jordan Pairs,''
Lecture Notes, Irvine, 1977

\[MWZ99 Magyar, P., J. Weyman, and A. Zelevinsky, {\it Multiple flag
varieties of finite
type}, Advances in Math. {\bf 141} (1999), 97--118

\[MWZ00 ---, {\it Symplectic multiple flag varieties of finite type},
J. Algebra {\bf 230:1}  (2000), 245--265

\[N\O{}04 Neeb, K-H., and B. \O{}rsted, {\it A topological Maslov
index for $3$-graded Lie groups},
Preprint TU Darmstadt {\bf 2366}, Oct. 2004

\[PR97 Popov, V., and G. R\"ohrle, {\it On the number of orbits of a
parabolic subgroup
on its unipotent radical}, in ``Algebraic Groups and Lie Groups. A volume of
papers in honour of the late R. W. Richardson,'' eds. Gus et al,
Cambridge Univ. Press,
Aust. Math. Soc. Lect. Ser. {\bf 9} (1997), 297--320

\[RRS92 Richardson, R., G. R\"ohrle, and R. Steinberg, {\it Parabolic
subgroups with
Abelian unipotent radical}, Invent. Math. {\bf 110} (1992), 649--671

\[Sa80 Satake, I., ``Algebraic Structures of Symmetric Domains," Publications
of the Math. Soc. of Japan {\bf 14}, Princeton Univ. Press, 1980

\[VR76 Vergne, M., and H.\ Rossi, {\it Analytic
continuation of the holomorphic
discrete series of a semisimple Lie group}, Acta Math. {\bf 136}(1976), 1--59

\[Vi86 Vinberg, E. B., {\it Complexity of the actions of reductive
groups}, Functional Anal. Appl.
{\bf 20} (1986), 1--13

\[Wh57 Whitney, H., {\it Elementary structure of real algebraic varieties},
Ann. of Math. (2) {\bf 66} (1957), 545--556

\[Wo71 Wolf, J. A., {\it Remark on Siegel domains of type III}, Proc.
Amer. Math. Soc. {\bf 30} (1971),
487--491

\[WK65 Wolf, J. A., and A. Kor\'anyi, {\it Realization of hermitean symmetric
spaces as generalized half planes}, Ann. of. Math. {\bf 81} (1965), 265--288

}

\references

\dlastpage
\bye